\documentclass[a4paper,10pt]{article}

\usepackage{amssymb}
\usepackage{graphicx}
\usepackage{graphics}

\def\yen{\hbox{iftdir\yoko\fi
\setbox0=\hbox{Y}Y\kern-.97\wd0\vbox{\hrule height.lex width.98\wd0
\kern.33ex\hrule height.lex width.98\wd0\kern.45ex}}}

\def\yen{{\setbox0=\hbox{Y}Y\kern-.97\wd0\vbox{hrule height.lex width.98%
\wd0\kern.33ex\hrule height.lex width.98\wd0\kern.45ex}}}

\def\np{\newpage}    
\begin{document}

\newtheorem{thm}{Theorem}[section]
\newtheorem{prop}{Lemma}

\def\g{\gamma}    
\def\ve{\varepsilon}   
\def\vp{\varphi}
\def\si{\mathrm{sin}} 
\def\co{\mathrm{cos}}
\def\i{\hookrightarrow}
\def\e{\hookrightarrow}
\def\l{\longrightarrow} 
\def\ttt{\longmapsto}
\def\f{\flushpar }
\def\nl{\newline }
\def\np{\newpage }
\def\x{\times } 
\def\te{^t \hskip-1mm }

\title{
Ribbon-moves of 2-links preserve \\
the $\mu$-invariant of 2-links 
}
\author{ Eiji Ogasa\\
 ogasa@ms.u-tokyo.ac.jp\\
Department of Physics, 
 University of Tokyo\\ 
 Hongo, Tokyo 113, JAPAN\\
}
\date{}
\maketitle

\noindent{\bf Abstract.}
We introduce ribbon-moves of 2-knots, which are operations to make 2-knots 
into new 2-knots by local operations in $B^4$. 
(We do not assume the new knots is not equivalent to the old ones.)

Let $L_1$ and $L_2$ be 2-links.  Then the following hold. 

(1) If $L_1$ is ribbon-move equivalent to $L_2$, 
then we have 

$$\mu(L_1)=\mu(L_2)$$.    

(2) 
Suppose that $L_1$ is ribbon-move equivalent to $L_2$. 
Let $W_i$ be arbitrary Seifert hypersurfaces for $L_i$. 
Then the torsion part of $H_1(W_1)\oplus H_1(W_2)$ is congruent to 
$G\oplus G$ for a finite abelian group $G$. 

(3) 
 Not all 2-knots are ribbon-move equivalent to the trivial 2-knot.   

(4) The inverse of (1) is not true. 

(5) The inverse of (2) is not true.

Let $L=(L_1,L_2)$ be a sublink of homology boundary link. 
Then we have: 
(i) $L$ is ribbon-move equivalent to a boundary link. 
(ii) $\mu(L)= \mu(L_1) + \mu(L_2)$.

We would point out the following facts 
by analogy of the discussions of finite type invariants of 1-knots 
although they are very easy observations. 
By the above result (1), we have: 
the $\mu$-invariant of 2-links is an order zero finite type invariant 
associated with ribbon-moves 
and 
there is a 2-knot whose $\mu$-invariant is not zero.
The mod 2 alinking number of $(S^2, T^2)$-links is 
an order one finite type invariant associated with the ribbon-moves
and there is an $(S^2, T^2)$-link whose mod 2 alinking number is not zero.

\newpage
\noindent{\large{\bf\S1. Introduction}}

In this paper we discuss ribbon-moves.  

An  {\it (oriented) (ordered) $m$-component 2-(dimensional) link}
 is a smooth, oriented submanifold $L=\{K_1,...,K_m\}$ of $S^4$, 
which is the ordered disjoint union of $m$ manifolds, each diffeomorphic  
to the $2$-sphere. 
 If $m=1$, then $L$ is called a {\it 2-knot}.  
We say that 2-links $L_1$ and $L_2$ are {\it equivalent} 
if there exists an orientation preserving diffeomorphism 
$f:$ $S^4$ $\rightarrow$ $S^4$ 
such that $f(L_1)$=$L_2$  and 
that $f | _{L_1}:$ $L_1$ $\rightarrow$ $L_2$ is 
an order and orientation preserving diffeomorphism.    
Let $id:S^4$ $\rightarrow$ $S^4$ be the identity. 
We say that 2-links $L_1$ and $L_2$ are {\it identical}  
if  $id(L_1)$=$L_2$  and 
that $id | _{L_1}:L_1$ $\rightarrow L_2$ is 
an order and orientation preserving diffeomorphism.    


We define ribbon-moves of 2-links.

\vskip3mm
\noindent
{\bf Definition 1.1.}
Let $L_1=(K_{1,1}...K_{1,m})$ 
and $L_2=(K_{2,1}...K_{2,m})$ 
be 2-knots in $S^4$. 
We say that $L_2$ is obtained from $L_1$ by one {\it ribbon-move } 
if there is a 4-ball $B$ of $S^4$ with the following properties.  

(1) 
$L_1-(B\cap L_1)=L_2-(B\cap L_2)$. 

 $K_{1,j}-(B\cap K_{1,j})=K_{2,j}-(B\cap K_{2,j})$ 

These diffeomorphism maps are orientation preserving.

(2) 
$B\cap L_1$ is drawn as in 
Figure 1.1.  
$B\cap L_2$ is drawn as in 
Figure 1.2.  

We regard $B$ as 
(a close 2-disc)$\times[0,1]\times\{t| -1\leq t\leq1\}$.
We put $B_t=$(a close 2-disc)$\times[0,1]\times\{t \}$.  
Then $B=\cup B_t$. 
In Figure 1.1 
and 1.2, 
we draw $B_{-0.5}, B_{0}, B_{0.5}$ $\subset B$. 
We draw $L_1$ and $L_2$ by the bold line. 
The fine line denotes $\partial B_t$. 
  
$B\cap L_1$ (resp. $B\cap L_2$) is diffeomorphic to 
$D^2\amalg (S^1\times [0,1])$.

$B\cap L_1$ has the following properties:  
$B_t\cap L_1$ is empty for $-1\leq t<0$ and $0.5<t\leq1$.
$B_0\cap L_1$ is diffeomorphic to 
$D^2\amalg(S^1\times [0,0.3])\amalg(S^1\times [0.7,1])$. 
$B_{0.5}\cap L_1$ is diffeomorphic to $(S^1\times [0.3,0.7])$. 
$B_t\cap L_1$ is diffeomorphic to $S^1\amalg S^1$ for $0<t<0.5$.

$B\cap L_2$ has the following properties:  
$B_t\cap  L_2$ is empty for $-1\leq t<-0.5$ and $0<t\leq1$.
$B_0\cap L_2$ is diffeomorphic to 
$D^2\amalg(S^1\times [0, 0.3])\amalg(S^1\times [0.7, 1])$. 
$B_{-0.5}\cap  L_2$ is diffeomorphic to $(S^1\times [0.3, 0.7])$. 
$B_t\cap  L_2$ is diffeomorphic to $S^1\amalg S^1$ for $-0.5<t<0$. 

We do not assume which the orientation of $B\cap L_1$ (resp. $B\cap L_2$ ) is. 

\hskip3cm Figure 1.1.   

\hskip3cm Figure 1.2. 

Suppose that $L_2$ is obtained from $L_1$ by one ribbon-move 
and that $L'_2$ is equivalent to $L_2$.   
Then we also say that $L'_2$ is obtained from $L_1$ 
by one {\it ribbon-move}.   
If $L_1$ is obtained from $L_2$ by one ribbon-move,  
then we also say that $L_2$ is obtained from $L_1$ by one {\it ribbon-move}.

\noindent
{\bf  Definition 1.2.}
2-knots $L_1$ and $L_2$ are said to be {\it ribbon-move equivalent} 
if there are 2-knots 
$L_1=\bar{L}_1, \bar{L}_2,...,\bar{L}_{p-1},\bar{L}_p=L_2$  
 ($p\in{\bf N}, p\geq2$) such that 
$\bar{L}_i$ is obtained from $\bar{L}_{i-1}$ $(1< i\leq p)$ 
by one ribbon-move.

In this paper we discuss  the following problems.

\vskip3mm
\noindent     {\bf Problem 1.3. }  
Let $L_1$ and $L_2$ be 2-links.
Consider a necessary (resp. sufficient, necessary and sufficient )
condition that $L_1$ and $L_2$ are ribbon-move equivalent.  
In particular, is there 
 a 2-knot which is not ribbon-move equivalent to the trivial 2-knot?
\vskip3mm

 \noindent
     {\bf Note.} 
(1)
Of course all $m$-component ribbon 2-links 
are ribbon-move equivalent to the trivial $m$-component link. 
See Appendix for the definition of ribbon 2-links.

(2)
By using \S4 of 
\cite{C1}, 
it is easy to prove that there is a nonribbon 2-knot 
which is ribbon-move equivalent to the trivial 2-knot.

Our motivation is as follows.  
We hope to investigate  `link space'
$E=\{f |  f:S^2\amalg...\amalg S^2\hookrightarrow S^4$ embeddings $\}$.  
In the case of 1-dimensional knots and links, 
we know that  
it is useful to investigate the space of immersions of circles 
in order to help investigate the space of embeddings. 
To discuss the space of immersions and that of embeddings 
is to discuss local moves (or knotting operations). 
In the case of 1-dimensional knots and links, 
we find many relations among `link space,' 
local moves, invariants of links, and QFT. 
(See 
\cite{J}    
\cite{Ko}   
\cite{V}    
\cite{W}    
etc.)
In 1-dimensional case, it is easy to find an unknotting operation. 
But high dimensional case, our first task is to define 
what kind of local moves we use. 
In this paper we discuss ribbon-moves as one of such moves.

This article is based on 
\cite{O2}. 
After 
\cite{O2},   
the author discusses relations between 
ribbon-moves of 2-knots and the Levie-Farber pairing 
and the Atiyah-Patodi-Singer-Casson-Gordon-Ruberman 
${\widetilde{\eta}}$-invariants of 2-knots (see 
\cite{O4} 
). 
In 
\cite{O3} 
the author discussed relations between 
local moves of $n$-knots and some invariants of $n$-knots.

\vskip3mm
\noindent{\large{\bf\S2. Main results}}

\vskip3mm
\noindent
{\bf Theorem 2.1 } {\it
Let $L_1$ and $L_2$ be 2-links in $S^4$. 
Suppose that $L_1$ is obtained from  $L_2$ by one ribbon-move.  
Then there are Seifert hypersurfaces 
$V_1$ for $L_1$ and $V_2$ for $L_2$   
such that 
($V_1$, $\sigma_1$) is spin preserving diffeomorphic to ($V_2$, $\sigma_2$), 
where  $\sigma_i$ is a spin structure induced from the unique one on $S^4$. 
} 
\vskip3mm

By using Theorem 2.1, we prove Theorem 2.2 and 2.3.   

\vskip3mm
\noindent
{\bf Theorem 2.2.} 
{\it
If 2-links $L$ and $L'$ are ribbon-move equivalent, then $\mu(L)$=$\mu(L')$. 
}
\vskip3mm

In \S 3 we define the $\mu$-invariant of 2-links. 

\vskip3mm
\noindent
{\bf Theorem 2.3.} 
{\it
Let $L_1$ and $L_2$ be 2-links in $S^4$. 
Suppose that $L_1$ are ribbon-move equivalent to $L_2$.   
Let $W_i$ be arbitrary Seifert hypersurfaces for $L_i$. 
Then the torsion part of $H_1(W_1)\oplus H_1(W_2)$ is congruent to 
$G\oplus G$ for a finite abelian group $G$. 
} 
\vskip3mm

By using Theorem 2.2 we prove Corollary 2.4.  
By using Theorem 2.3 we also prove Corollary 2.4.    

\vskip3mm
\noindent
{\bf Corollary 2.4. }   
{\it 
Not all 2-knots are ribbon-move equivalent to the trivial 2-knot.   
}
\vskip3mm

By using Theorem 2.3 we prove Corollary 2.5.    

\vskip3mm
\noindent
{\bf Corollary 2.5. }   
{\it
There is a 2-knot $K$ such that 
$\mu(K)=0$ and that   
$K$ is not ribbon-move equivalent to the trivial 2-knot.   
}
\vskip3mm

By using Theorem 2.2,  we prove Corollary 2.6.    

\vskip3mm
\noindent
{\bf Corollary 2.6. }   
{\it
The inverse of Theorem 2.3 is not true.   
}
\vskip3mm

In \S3-7 we prove the above results.

In \S8 we prove that: 
Let $L=(L_1,L_2)$ be a sublink of homology boundary link. 
Then the following hold. 
(1) $L$ is ribbon-move equivalent to a boundary link. 
(2) $\mu(L)= \mu(L_1) + \mu(L_2)$. 

In \S9 we would point out the following facts 
by analogy of the discussions of 
finite type invariants of 1-knots although 
they are very easy observations. 
By Theorem 2.2, we have: 
the $\mu$-invariant of 2-links is an order zero finite type invariant 
associated with ribbon-moves 
and 
there is a 2-knot whose $\mu$-invariant is not zero.
The mod 2 alinking number of $(S^2, T^2)$-links is 
an order one finite type invariant associated with the ribbon-moves
and there is an $(S^2, T^2)$-link whose mod 2 alinking number is not zero.

\vskip3mm
\noindent{\large{\bf\S3. 
The $\mu$-invariant of 2-links }}  

See \S IV of 
\cite{Ki}   
 for the spin structures and 
 the $\mu$-invariant of closed spin 3-manifolds.

\vskip3mm
\noindent{\bf Definition.}
Let $L=(K_1,...,K_m)$ be a 2-link.
Let $V$ be a Seifert hypersurface for $L$. 
Note that $V$ is oriented so that 
the orientation is compatible with that on $L$.
A spin structure $\sigma$ on $V$ is induced 
from the unique spin structure on $S^4$. 
Attach $m$ 3-dimensional 3-handles to $V$ along each component of 
the boundary. 
Then we obtain the closed oriented 3-manifold $\hat V$.   
The spin structure $\sigma$ extends over $\hat V$ uniquely.
Call it $\hat\sigma$.
We define 
the {\it $\mu$-invariant $\mu(L)$ of the 2-link $L$}     
to be the $\mu$-invariant $\mu((\hat V, \hat\sigma))$ 
$\in {\bf Z_{16}}$ 
of the closed spin 3-manifold $(\hat V, \hat\sigma)$. 
\vskip3mm

\vskip3mm
\noindent
{\bf Claim.}   {\it 
Under the above conditions 
$\mu(L)$ is independent of the choice of $V$.      
}
\vskip3mm

\noindent{\bf Proof.} 
P.580 of 
\cite{R}  
proved the above Claim when $L$ is a knot.

\cite{Kw}   
 says:

\vskip3mm
\noindent
{\bf Fact 3.1. }   
{\it
(\cite{Kw}) 
Let $V$ and $V'$ be Seifert hypersurfaces for $L$. 
Then we have: 
there are Seifert hypersurfaces 
$V=V_1$,$V_2$,...,$V_{p-1}$,$V_p$ for $L$ 
with the following properties.  

(1)The embedding map of $V_p$ is isotopic to that of $V'$, 
where we do not fix the boundary of the image. 
(Note. $[V\cup V']$ is not zero in general in  $H_3(S^4-L;{\bf Z})$. 
But we can set $[V_1\cup V_p]=0$ $\in H_3(S^4-L; {\bf Z})$. )

(2) 
For $V_i$ and $V_{i+1}$ ($i=1,...,p-1$),  
there is a compact oriented 4-manifold $W_i$ embedded in $S^4$ 
which has a handle decomposition 

$W_i=
(V_i\times[0,1])\cup$ (one $q$-handle) $\cup( V_{i+1}\times[0,1])$ 
($q\in\{1,2,3\}$). 
}
\vskip3mm

We give $W_i$ a spin structure induced from the unique one on $S^4$.

The following two spin structures on $V_1$ coincide one another. 
Call it $\sigma_1$. 

(i) The spin structure induced from the unique one on $S^4$ 

(ii) The spin structure induced from the one on $W_1$.

The following two spin structures on $V_p$ coincide one another. 
 Call it $\sigma_{p}$. 

(i)  The spin structure induced from the unique one on $S^4$ 

(ii)  The spin structure induced from the one on $W_{p-1}$.

The following three spin structures on $V_i$ 
coincide each other ($i=2,...,p-1$).  
Call it $\sigma_i$.

(i) The spin structure induced from the unique one on $S^4$.

(ii) The spin structure induced from the one on $W_i$.

(iii) The spin structure induced from the one on $W_{i+1}$.

The 3-dimensional closed oriented spin 3-manifolds  
$(\hat V_i,\hat\sigma_i)$ are defined from $(V_i,\sigma_i)$ as 
in the above Definition ($i=1,...,p$).   
(See \S IV of 
\cite{Ki} 
for the way to induce spin structures on manifolds 
from those on others. 
)

Let $x,y$ be arbitrary elements of $H_2(W_i;{\bf Z})/{\mathrm{Tor}}$. 
Let $x\cdot y$ be the intersection product.  

We prove: $x\cdot y=0.$ 

There is an oriented closed surface $F$ embedded in $W_i$ 
which represents $x$. 
Since $F$ is embedded in $S^4$, $[F]\cdot[F]=0$. 
Hence $x\cdot x=0$ for any element 
$x\in H_2(W_i;{\bf Z})/{\mathrm{Tor}}$.
Hence $x\cdot y=0$ for arbitrary elements  
$x, y\in H_2(W_i;{\bf Z})/{\mathrm{Tor}}$.
Hence the signature of the intersection form  

$H_2(W_i;{\bf Z})/{\mathrm{Tor}}\times H_2(W_i;{\bf Z})/{\mathrm{Tor}}
\to{\bf Z}$  \quad $(x,y)\mapsto x\cdot y$

is the zero map. 
Hence $\sigma(W_i)=0$. 

Therefore 
$\mu((\hat V_i,\hat\sigma_i))-\mu(-(\hat{V}_{i+1},\hat{\sigma}_{i+1}))$=
$\mu((V_i,\sigma_i)\cup (-(V_{i+1},\sigma_{i+1}))$=

\{mod 16 $\sigma(W_i)\}=0$. 
Hence 
$\mu((\hat V_i,\hat\sigma_i))=\mu(\hat{V}_{i+1},\hat{\sigma}_{i+1})$. 
$(i=1,...p-1.)$

Therefore 
$\mu((\hat V_1,\hat\sigma_1))$=
$\mu((\hat V_2,\hat\sigma_2))$=...=
$\mu((\hat V_{p-1},\hat\sigma_{p-1}))$=
$\mu(\hat{V}_p,\hat{\sigma}_p)$.

This completes the proof.

\vskip3mm
\noindent{\large{\bf\S4.  Proof of Theorem 2.1 }}    

In order to prove Theorem 2.1,   
we introduce (1,2)-pass-moves of 2-links.  

\noindent{\bf Definition 4.1.}  
Let $L_1=(K_{1,1}...K_{1,m})$ 
and $L_2=(K_{2,1}...K_{2,m})$ 
be 2-knots in $S^4$. 
We say that $L_2$ is obtained from $L_1$ by one {\it (1,2)-pass-move } 
if there is a 4-ball $B$ $\subset S^4$ with the following properties.  
We draw $B$ as in 
Definition 1.1. 

(1) 
$L_1-(B\cap L_1)$=$L_2-(B\cap L_2)$. 

 $K_{1,j}-(B\cap K_{1,j})=K_{2,j}-(B\cap K_{2,j})$

These diffeomorphism maps are orientation preserving. 

(2)
$B\cap L_1$  is drawn as in 
Figure 4.1. 
$B\cap L_2$ is drawn as in 
Figure 4.2.  

\hskip3cm Figure 4.1.   

\hskip3cm Figure 4.2.   

The orientation of the two discs in the 
Figure 4.1 (resp. Figure 4.2)   
is compatible with 
the orientation which is naturally determined 
by the $(x,y)$-arrows in the Figure. 
We do not assume which the orientations of the annuli in the Figures are.

Suppose that $L_2$ is obtained from $L_1$ by one (1,2)-pass-move 
and that $L'_2$ is equivalent to $L_2$.   
Then we also say that $L'_2$ is obtained from $L_1$ by one 
{\it (1,2)-pass-move }.

If $L_1$ is obtained from $L_2$ by one (1,2)-pass-move,  then 
we also say that $L_2$ is obtained from $L_1$ by one {\it (1,2)-pass-move }.  

2-knots $L_1$ and $L_2$ are said to be {\it (1,2)-pass-move equivalent} 
if there are 2-knots 
$L_1=\bar{L}_1, \bar{L}_2,...,\bar{L}_{p-1},\bar{L}_p=L_2$  
$( p\in {\bf N}, p\geq2 )$ such that 
$\bar{L}_i$ is obtained from $\bar{L}_{i-1}$ $(1< i\leq p)$ 
by one (1,2)-pass-move.

\vskip3mm
\noindent
{\bf Proposition 4.2.}   
 {\it
Let $L$ and $L'$ be 2-links.  
Then the  following conditions (1) and (2) are equivalent. 

(1)
$L$ is  (1,2)-pass-move equivalent to $L'$.

(2)
$L$ is  ribbon-move equivalent to $L'$. 
}
\vskip3mm

It is obvious that 
Proposition 4.2 
follows from 
Proposition 4.3. 

\vskip3mm
\noindent
{\bf Proposition 4.3.}    
 {\it
Let $L$ and $L'$ be 2-links.  Then the following hold.

(1)
If $L$ is  obtained from $L'$ by one ribbon-move, 
then $L'$ is  obtained from  $L$ by one (1,2)-pass-move.

(2)
If $L$ is  obtained from $L'$ by one (1,2)-pass-move, 
then $L'$ is  obtained from  $L$ by two ribbon-move.
}
\vskip3mm

Proposition 4.3.(2) 
is obvious. 

Proposition 4.3.(1)  
 follows from 
Proposition 4.4  
because: 
The pair of a manifold and a submanifold, 
$($ the 4-ball, $($the 2-link$)\cap($the 4-ball$)$$)$, 
in Figure 4.1 is included in 
the pair 
$($ the 4-ball, $($the 2-link$)\cap($the 4-ball$)$$)$ 
in  Figure 4.4.

\vskip3mm
\noindent
{\bf Proposition 4.4.}   
{\it
Let $L_1=(K_{1,1},...,K_{1,m})$ 
and $L_2=(K_{2,1},...,K_{2,m})$ be 2-links in $S^4$. 
Then the following two conditions $(I)$ and $(II)$ are equivalent. 

$(I)$ $L_1$ is equivalent to $L_2$. 

$(II)$ There is 
a 4-ball $B$ $\subset S^4$ with the following properties.  
We draw  $B$ as in 
Definition 1.1. 

(1) 
$L_1-(B\cap L_1)=L_2-(B\cap L_2)$. 

$K_{1,i}-(B\cap K_{1,i})=K_{2,i}-(B\cap K_{2,i})$ for each $i$. 

These diffeomorphism maps are orientation preserving.

(2)   
$B\cap L_1$ is drawn as in Figure 4.3.  
$B\cap L_2$ is drawn as in Figure 4.4.  

\hskip3cm Figure 4.3.    

\hskip3cm Figure 4.4.    

The orientation of $B\cap L_2$ is compatible with 
the orientation which is naturally determined 
by the $(x,y)$-arrows in the Figure 4.4. 
}
\vskip3mm

It is obvious that  
Proposition 4.4  
follows from 
Proposition 4.5. 

\vskip3mm
\noindent
{\bf  Proposition 4.5. }  
 {\it
Let $L_1=(K_{1,1},...,K_{1,m})$ 
and $L_2=(K_{2,1},...,K_{2,m})$ be 2-links in $S^4$. 
Then the following two conditions $(I)$ and $(II)$ are equivalent. 

$(I)$ $L_1$ is equivalent to $L_2$. 

$(II)$ There is a 4-ball $B \subset S^4$ with the following properties.  
We draw $B$ as in 
Definition 1.1. 
 
(1)
$L_1-(B\cap L_1)=L_2-(B\cap L_2)$. 

$K_{1,i}-(B\cap K_{1,i})$=$K_{2,i}-(B\cap K_{2,i})$ for each $i$.

These diffeomorphism maps are orientation preserving.

(2)
$B\cap L_1$ is drawn as in Figure 4.5.  
$B\cap L_2$ is drawn as in Figure 4.6.  

\hskip3cm Figure 4.5.  

\hskip3cm Figure 4.6.  

We do not assume which the orientation of $B\cap L_1$ (resp. $B\cap L_2$ ) is. 
}
\vskip3mm

Proposition 4.5  
follows from 
Proposition 4.6 
because: 
The pair of a manifold and a submanifold, 
$($ the 4-ball, $($the 2-link$)\cap($the 4-ball$)$$)$, 
in Figure 4.5 (resp. Figure 4.6) 
is made from 
the pair 
$($ the 4-ball, $($the 2-link$)\cap($the 4-ball$)$$)$ 
in  Figure 4.7 (resp. Figure 4.8) 
by a rotation through $90^\circ$ around an appropriate plane in the 4-ball 
in Figure 4.5 (resp. Figure 4.6) 
and by isotopy. 
 

\vskip3mm
\noindent
{\bf  Proposition 4.6.} 
 {\it
Let $L_1=(K_{1,1},...,K_{1,m})$ 
and $L_2=(K_{2,1},...,K_{2,m})$ be 2-links in $S^4$. 
Then the following two conditions $(I)$ and $(II)$ are equivalent. 

$(I)$ $L_1$ is equivalent to $L_2$. 

$(II)$ There is 
a 4-ball $B \subset S^4$ with the following properties.  
We draw $B$ as in 
Definition 1.1. 
 
(1)
$L_1-(B\cap L_1)=L_2-(B\cap L_2)$. 

$K_{1,i}-(B\cap K_{1,i})$=$K_{2,i}-(B\cap K_{2,i})$ for each $i$. 

These diffeomorphism maps are orientation preserving.

(2)
$B\cap L_1$ is drawn as in Figure 4.7.  
$B\cap L_2$ is drawn as in Figure 4.8.  

\hskip3cm Figure 4.7.  

\hskip3cm Figure 4.8.  

We do not assume which the orientation of $B\cap L_1$ (resp. $B\cap L_2$ ) is. 
}
\vskip3mm

\noindent{\bf Proof of Proposition 4.6.} 
We obtain $L_2$ from $L_1$ by an explicit isotopy,  
 Figure 4.7 $\to$ Figure 4.9 $\to$ Figure 4.8.   
Note that the following 
Proposition 4.7 
holds by an explicit isotopy. 
This  completes the proof of Proposition 4.2-4.6. 

\hskip3cm Figure 4.9.  

\vskip3mm
\noindent
{\bf  Proposition 4.7.} 
 {\it
Let $L_1=(K_{1,1},...,K_{1,m})$ 
and $L_2=(K_{2,1},...,K_{2,m})$ be 2-links in $S^4$. 
Then the following two conditions $(I)$ and $(II)$ are equivalent. 

$(I)$ $L_1$ is equivalent to $L_2$. 

$(II)$ There is a 4-ball $B \subset S^4$ with the following properties.  
We draw $B$ as in 
Definition 1.1. 
 
(1)
$L_1-(B\cap L_1)=L_2-(B\cap L_2)$. 

$K_{1,i}-(B\cap K_{1,i})$=$K_{2,i}-(B\cap K_{2,i})$ for each $i$. 

These diffeomorphism maps are orientation preserving.

(2)
$B\cap L_1$ is drawn as in Figure 4.10.  
$B\cap L_2$ is drawn as in Figure 4.11.  

\hskip3cm Figure 4.10.  

\hskip3cm Figure 4.11.  

We do not assume which the orientation of $B\cap L_1$ (resp. $B\cap L_2$) is. 
}
\vskip3mm

\noindent{\bf Note.}    
Regard the operation,  

`t=0 of Figure 4.7 $\to$ t=0 of Figure 4.8 $\to$ t=0 of Figure 4.9,' 

\noindent 
as an isotopy of (a part of) 1-knot. 
Then this operation is essentially same as 
the operation in the figure in the proof of 
Lemma 5.5 
of 
\cite{Kf}. 

\noindent{\bf Proof of Theorem 2.1. }   
By Proposition 4.3.(1),   
$L_1$ is obtained from $L_2$ by one (1,2)-pass-move in a 4-ball $B$.

\vskip3mm
\noindent
{\bf Claim 4.8.} {\it 
There are Seifert hypersurfaces $V_1$ for $K_1$ and $V_2$ for $K_2$ 
such that: 


(1) 
$V_1-(B\cap V_1)=V_2-(B\cap V_2)$.

These diffeomorphism maps are orientation preserving.

(2) 
$B\cap V_1$ is drawn as in Figure 4.12. 
$B\cap V_2$ is drawn as in Figure 4.13. 
}
\vskip3mm

\noindent
{\bf Note.}  
We draw $B$ as in Definition 1.1. 
We draw $V_1$ and $V_2$ by the bold line. 
The fine line means $\partial B$. 
  
$B\cap V_1$ (resp. $B\cap V_2$) is diffeomorphic to 
$(D^2\times [2,3])\amalg (D^2\times [0,1])$. 
We can regard $(D^2\times [0,1])$ as a 3-dimensional 1-handle 
which is attached to $\partial B$. 
We can regard $(D^2\times [2,3])$ as a 3-dimensional 2-handle 
which is attached to $\partial B$.

$B\cap V_1$ has the following properties:  
$B_t\cap V_1$ is empty for $-1\leq t<0$ and $0.5<t\leq1$.
$B_0\cap V_1$ is diffeomorphic to 
$(D^2\times [2,3])\amalg(D^2\times [0,0.3])\amalg(D^2\times [0.7,1])$. 
$B_{0.5}\cap K_1$ is diffeomorphic to $(D^2\times [0.3,0.7])$. 
$B_t\cap V_1$ is diffeomorphic to $D^2\amalg D^2$ for $0<t<0.5$.

$B\cap V_2$ has the following properties: . 
$B_t\cap  V_2$ is empty for $-1\leq t<-0.5$ and $0<t\leq1$.
$B_0\cap V_2$ is diffeomorphic to 
$(D^2\times[2,3])\amalg(D^2\times [0, 0.3])\amalg(D^2\times [0.7, 1])$. 
$B_{-0.5}\cap  V_2$ is diffeomorphic to $(D^2 \times [0.3, 0.7])$. 
$B_t\cap V_2$ is diffeomorphic to $D^2\amalg D^2$ for $-0.5<t<0$.

\hskip3cm  Figure 4.12   

\hskip3cm  Figure 4.13. 

\noindent
{\bf Proof of Claim. }  
Put  $P=($the 3-manifolds in Figure 4.12$)\cap(\partial B)$. 
Note $P=($the 3-manifolds in Figure 4.13$)\cap(\partial B)$. 
Put $Q=L_1\cap(S^4-Int B^4)$. 
Note $Q=L_2\cap(S^4-Int B^4)$. 
By applying the following Proposition to ($P\cup Q$) 
and $(S^4-Int B^4)$, Claim 4.8 holds.

The following Proposition is proved by using the obstruction theory. 
We give a proof although it is folklore.  

\vskip3mm
\noindent
{\bf Proposition.} 
{\it 
Let $X$ be an oriented compact $(m+2)$-dimensional manifold. 
Let $\partial X\neq\phi$. 
Let $M$ be an oriented closed $m$-dimensional manifold 
which is embedded in $X$. 
Let $M\cap\partial X\neq\phi$. 
Let $[M]=0\in H_m(X; {\bf Z})$. 
Then there is an oriented compact $(m+1)$-dimensional manifold $P$ 
such that $P$ is embedded in $X$ and that $\partial P=X$. 
}
\vskip3mm

\noindent{\bf Proof.} 
Let $\nu$ be the normal bundle of $M$ in $X$.  
By Theorem 2 in P.49 of 
\cite{Ki}   
$\nu$ is a product bundle. 
By using $\nu$ and the collar neighborhood of $\partial X$ in $X$, 
we can take a compact oriented $(m+2)$-manifold $N\subset X$ 
with the following properties. 

(1)
 $N\cong M\times D^2$. (Hence $\partial N=M\times S^1$.)

(2) 
$N\cap\partial X$
$=(\partial N)\cap(\partial X)$
$=M\cap\partial X$. 
(Hence  (Int$N)\cap\partial X=\phi$. )

Take $X-($Int$N)$. (Note $X-($Int$N)\supset \partial X$.)
There is a cell decomposition: 

$X-({\mathrm{Int}}N)$

$=(\partial N)\cup(\partial X)\cup$
(1-cells $e^1)\cup$  
(2-cells $e^2)\cup$ 
(3-cells $e^3)\cup$ 
(one 4-cell $e^4)$.

We can suppose that this decomposition has only one 0-cell $e^0$ 
which is in $(\partial N)\cap(\partial X)$.


There is a continuous map 
$s_0:(\partial N)\cup(\partial X)\to S^1$ 
with the following properties, 
where $p$ is a point in $S^1$.  

(1)  $s_0(\partial X)=p$. 
    (Hence $s_0((\partial N)\cap(\partial X))=p$ and 
    $s_0(e^0)=p$. )

(2) $s_0\vert_{\partial N}: M\times S^1\to S^1$ is a projection map 
     $(x,y)\mapsto y$.

Let $S_F^1$ be a fiber of the $S^1$-fiber bundle $\partial N=M\times S^1$. 
Since $[M]=0\in H_m(X; {\bf Z})$, 
$[S_F^1]$ generates ${\bf Z}\subset$ $H_1(X-$Int $N, \partial X;{\bf Z})$. 
(We can prove as in the proof of Theorem 3 in P.50 of 
\cite{Ki}   
) 

Let $f:H_1(X-\mathrm{Int}N, \partial X;{\bf Z})
\to H_1(X-\mathrm{Int}N, \partial X;{\bf Z})/$Tor 
be the natural projection map. 
Let $\{f([S^1_F]), u_1,...,u_k\}$ be a set of basis of 

\noindent$H_1(X-\mathrm{Int}N, \partial X;{\bf Z})/$Tor. 
Take a continuous map 

$s_1: (\partial N)\cup(\partial X)\cup$(1-cells $e^1)\to S^1$  

with the following properties. 

(1) $s_1\vert_
{(\partial N)\cup(\partial X)}=s_0$  

(2) 
$s_1\vert_{e^0\cup e^1}:e^0\cup e^1\to S^1$ 
satisfies the following condition: 
If $f([e^0\cup e^1])=$
\newline
$n_0\cdot f([S^1_f])+\Sigma_{j=1}^{k}n_j\cdot u_j$
$\in H_1(X-\mathrm{Int}N, \partial X;{\bf Z})/$Tor
($n_*\in{\bf Z}$), 
\newline  
then deg$(s_1\vert_{e^0\cup e^1})=n_0$.

Note that, if a circle $C$ is nul-homologous in 
$(\partial N)\cup(\partial X)\cup$(1-cells $e^1)$, 
then deg$(s_1\vert_C)=0$. 

\noindent
{\bf Claim.}  
{\it 
There is a continuous map 

$s_2:(\partial N)\cup(\partial X)\cup$(1-cells $e^1)\cup$(2-cells $e^2)
\to S^1$  

such that 
$s_2\vert_{(\partial N)\cup(\partial X)\cup({\mathrm{1-cells}}\quad e^1)}$
=$s_1$.
}

\noindent
{\bf Proof.} 
It is trivial that $[\partial e^2]=0$
$\in H_1((\partial N)\cup(\partial X)\cup$(1-cells $e^1);{\bf Z})$. 
Hence deg$(s_1\vert_{\partial e^2})=0.$
Hence $s_1\vert_{\partial e^2}$ extends to $e^2$. 
Hence the above Claim holds.

The map $s_2$ extends to a continuous map $s:X-({\mathrm{Int}}N)\to S^1$ 
since $\pi_l(S^1)=0 (l\geq2)$. 
We can suppose $s$ is a smooth map.

Let $q\neq p$. Let $q$ be a regular value. 
Hence $s^{-1}(q)$ be an oriented compact manifold. 
$\partial\{s^{-1}(q)\}\subset\{(\partial N)\cup\partial X\}$.   
Since $q\neq p$, $s^{-1}(q)\cap\partial X=\phi.$  
Hence $\partial\{s^{-1}(q)\}\subset\partial N$.
Furthermore we have 
$s^{-1}(q)\cap \partial N=\partial \{s^{-1}(q)\} 
=M\times \{r\}$, where $r$ is a point in $S^1$. 
By using $N$ and $s^{-1}(q)$, Proposition holds.

By Claim 4.8, there is a smooth transverse immersion 
$F:V\times [1,2]\to S^4$ 
 such that 
 $F\vert_{V\times \{1\}}(V\times \{1\})=V_1$
 and 
  $F\vert_{V\times \{2\}}(V\times \{2\})=V_2$. 
Give a spin structure $\alpha$ on $V\times [1,2]$ 
by using $F$. 
Then the following two spin structures on $V_i$ coincide one another. 
Call it $\tau_i$.  

(i)  the spin structure induced from the unique spin structure $S^4$ 

(ii) the spin structure induced from 
 $\alpha$ on $V\times [1,2]$. 

By using $F$, 
it holds that $V_1$ and $V_2$ are spin preserving diffeomorphism. 
This completes the proof of 
Theorem 2.1. 

\vskip3mm
\noindent{\large{\bf\S5.  Proof of Theorem 2.2 }}  

By Proposition 4.2,   
$L$ and $L'$ are (1,2)-pass-move equivalent. 
Take 2-links 
$L=\bar{L}_1,\bar{L}_2,...,\bar{L}_{p-1},\bar{L}_p=L'$ 
as in 
Definition 4.1.  
Obviously, it suffices to prove that 
$\mu(\bar{L}_i)$=$\mu(\bar{L}_{i+1})$ for each $i$ ($1\leq i<p$). 
By Theorem 2.1   
we have: There are Seifert hypersurfaces,
 $V_{i,i+1}$ for $\bar{L}_i$ and $V_{i+1,i}$ for $\bar{L}_{i+1}$,  
such that $V_{i,i+1}$ and $V_{i+1,i}$ are  spin preserving diffeomorphism. 
 Hence $\mu(\bar{L}_i)$=$\mu(\bar{L}_{i+1})$.

\vskip3mm
\noindent{\large{\bf\S6.  The proof of Theorem 2.3  }}  

The following Fact 6.1   
is an elementary fact.

\vskip3mm
\noindent
{\bf  Fact6.1.} {\bf(Known)}  
 {\it
Let $A, B, C, X$ and $Y$ be a finite abelian group. 
Suppose $A\oplus B\cong X\oplus X$ and
$B\oplus C\cong Y\oplus Y$.  
Then  
$A\oplus C\cong P\oplus P$ for a finite abelian group $P$.   
}
\vskip3mm

It is obvious that Theorem 2.3 follows from
Theorem 2.1, Fact 6.1,    
and Proposition 6.2.  

\vskip3mm
\noindent
{\bf Proposition 6.2.}   
{\it
Let $V$ and $V'$ be Seifert hypersurfaces for  a 2-link $L$.  
Then the torsion part of 
$H_1(V;{\bf Z})\oplus H_1(V';{\bf Z})$ is congruent to $G\oplus G$  
for a finite group  $G$.
}
\vskip3mm

\noindent{\bf Proof.}  
Take $V_1,...,V_p$ and $W_1,...,W_p$ as in 
Fact 3.1  
and its proof. 
By using the Meyer-Vietoris sequence, we have 
Tor $H_1(\partial W_i;{\bf Z})\cong$ 
Tor $\{H_1(V_i;{\bf Z})\oplus H_1(V_{i+1};{\bf Z})\}$ 
($i=1,...,p-1$).
The manifold 
$\partial W_i$ is a closed oriented 3-manifold embedded in $S^4$. 
Hence 

\hskip1cm
Tor $H_1(\partial W_i;{\bf Z})\cong G_i\oplus G_i$ ------$(*)$

\noindent
for a finite abelian group $G_i$. 
 ( See e.g.
\cite{GL} 
\cite{H}.  
We give a proof in the following paragraph.)
Hence 
Tor $\{H_1(V_i;{\bf Z})\oplus H_1(V_{i+1};{\bf Z})\}$ 
$\cong G_i\oplus G_i$.

We give a proof for the above congruence $(*)$: 
By using the Meyer-Vietoris sequence 
$H_i(\partial W_i;{\bf Z})\to 
H_i(W_i;{\bf Z})\oplus H_i(\overline{S^4-W_i};{\bf Z})\to H_i(S^4;{\bf Z})$, 
Tor $H_1(\partial W_i;{\bf Z})\cong$ 
Tor $\{H_1( W_i;{\bf Z})\oplus H_1(\overline{S^4-W_i};{\bf Z})\}$. 
By using the Meyer-Vietoris sequence 
$H_i( W_i;{\bf Z})\to H_i(S^4;{\bf Z})\to 
H_i(S^4, W_i;{\bf Z})$,   
$H_1( W_i;{\bf Z})\cong H_2(S^4, W_i;{\bf Z}).$     
By the excision, $H_2(S^4, W_i;{\bf Z})\cong 
H_2(\overline{S^4-W_i}, \partial W_i;{\bf Z}).$     
By the Poincar\'e duality, 
$H_2(\overline{S^4-W_i}, \partial W_i;{\bf Z})\cong
H^2(\overline{S^4-W_i};{\bf Z}).$      
By the universal coefficient theorem, 
Tor$H_1(\overline{S^4-W_i};{\bf Z})\cong$
Tor$H^2(\overline{S^4-W_i};{\bf Z}).$      
Hence Tor$H_1(\overline{S^4-W_i};{\bf Z})\cong$ Tor$H_1(W_i;{\bf Z})$. 
Hence 
Tor$H_1(\partial W_i;{\bf Z})\cong$ 
Tor$H_1(\overline{S^4-W_i};{\bf Z})\oplus$ 
Tor$H_1(\overline{S^4-W_i};{\bf Z})$. 
Hence the congruence $(*)$ holds.

By Fact 6.1,     
Tor $\{H_1(V_1;{\bf Z})\oplus H_1(V_p;{\bf Z})\}$ 
=$G\oplus G$ 
for a finite abelian group $G$. 
Hence 
Tor $\{H_1(V;{\bf Z})\oplus H_1(V';{\bf Z})\}$ 
=$G\oplus G$.

\vskip3mm
\noindent{\large{\bf
\S7. The proof of Corollary 2.4, 2.5 and 2.6 }}  

Let $K$ be the 2-twist spun knot of a 1-knot $A$.
Let $M$ be the 2-fold branched cyclic covering space of $S^3$ along $A$. 
By 
\cite{Z}, 
$\overline{M-B^3}$ is a Seifert hypersurface for $K$. 
Let $S$ be a Seifert matrix of $K$. 
By Lemma 12.1, Theorem 12.2, and Theorem 12.6 
in Chapter XII of 
\cite{Kf},  
there is a compact oriented 4-manifold $X$ 
with the following properties. 
(i) $M=\partial X$. 
(ii) $H_1(X;{\bf Z})\cong0$
(iii) $H_3(X;{\bf Z})\cong0$
(iv) The intersection form $H_2(X;{\bf Z})\times H_2(X;{\bf Z})\to{\bf Z}$
  is represented by $S+^tS$. 
  (Note: 
  By using the Poincar\'e duality, the universal coefficient theorem,  
  and the above conditions (ii) (iii), 
   it holds that $H_2(X;{\bf Z})$ is torsion free.)

By the above fact (iv),  the intersection form is even. 
By this fact, the above (ii), and P.27 of
\cite{Ki},  
it holds that  $X$ is a spin manifold. 
Hence, for a spin structure $\alpha$ on $M$, 
$\mu(M,\alpha)=$ mod  16  $\sigma(S+^tS)$. 
 (Note that there is a spin 3-manifold 
whose spin structure is more than one.)

(1) 
Let $A$ be the trefoil knot. Let $S$ be 
$
\left(
\begin{array}{cc}
1&1\\
0&1\\
\end{array}
\right) 
$. 
Then the intersection form of $X$ is 
represented by 
$
\left(
\begin{array}{cc}
2&1\\
1&2\\
\end{array}
\right) 
$.

Hence we have: 

(1.1) $H_1(\overline{M-B^3};{\bf Z})\cong {\bf Z_3}$.

Hence $H_1(\overline{M-B^3};{\bf Z_2})\cong 0$. 
Hence $\overline{M-B^3}$ has only one spin structure. 
Hence $\mu(K)=\mu(M)=$ 
mod 16 
$ ( \sigma
\left(
\begin{array}{cc}
2&1\\
1&2\\
\end{array}
\right) 
)$.   
  Hence we have:  

(1.2)$\mu(K)=2$.

(2) Let $A$ be the figure eight knot. 
Let $S$ be 
$
\left(
\begin{array}{cc}
1&1\\
0&-1\\
\end{array}
\right) 
$. 
Then the intersection form of $X$ is 
represented by 
$
\left(
\begin{array}{cc}
2&1\\
1&-2\\
\end{array}
\right) 
$.

Then we have: 

(2.1)$H_1(\overline{M-B^3};{\bf Z})\cong {\bf Z_5}$. 

Hence $H_1(\overline{M-B^3};{\bf Z_2})\cong 0$.   
Hence $M$ has only one spin structure. 
Hence $\mu(K)=\mu(M)=$ 
mod 16 
$ ( \sigma
\left(
\begin{array}{cc}
2&1\\
1&-2\\
\end{array}
\right) 
)$. 
 Hence we have:

(2.2)  $\mu(K)=0$.

(3) Let $K$ be the 5-twist spun knot of the trefoil knot. 
Let $M$ be the Poincar\'e homology sphere.  
Then we have:  

(3.1) 
There is a Seifert hypersurface for $K$ which is diffeomorphic to 
$\overline{M-B^3}$. 
(See \S65 of 
\cite{ST}.  
)

(3.2) $\mu(K)=\mu(M)=8$. 
(See e.g. 
P.15 and P.67 of 
\cite{Ki}.  
)

The above (1.2) and Theorem 2.2 imply  Corollary 2.4.    

The above (1.1) (or (2.1)) and Theorem 2.3 imply  Corollary 2.4. 

The above (2.1), (2.2) and Theorem 2.3 imply Corollary 2.5.    

The above (3.1), (3.2) and Theorem 2.2 imply Corollary 2.6. 

\vskip3mm
\noindent{\large
{\bf\S8. Any SHB link is ribbon-move equivalent to a boundary link }}

See P.640 of 
\cite{C2}  
 and P.536 of 
\cite{CO}  
etc. 
for sublinks of homology boundary links 
( i.e. SHB links ), 
homology boundary links and boundary links.

\vskip3mm
\noindent
{\bf Theorem 8.1.}   
 {\it
Let $L=(K_1,K_2)$ be a 2-link. 
Let $L$ be a sublink of a homology boundary link.
Then $L$ is ribbon-move equivalent to a boundary link. 
}
\vskip3mm


To prove Theorem 8.1, we need lemmas. 
By the definition of SHB links (in P.536 of 
\cite{CO}   
) 
the following holds.  

\vskip3mm
\noindent
{\bf Lemma 8.1.1.}    
{\it
Let $L=(K_1,K_2)$ be the 2-link in 
Theorem 8.1.   
There is a connected Seifert hypersurface $V_i$ for $K_i$  ( $i=1,2$ )
such that 
$V_1\cap V_2$ is diffeomorphic to a disjoint union of 2-spheres 
$S^2_1$,...,$S^2_\nu$.
}
\vskip3mm

We prove: 
\vskip3mm
\noindent
{\bf Lemma 8.1.2.}   
{\it
Let $L=(K_1, K_2)$ be the 2-link in 
Theorem 8.1.   
Then there is a 2-link $L'=(K'_1,K'_2)$ which is equivalent to $L$ 
satisfying the following condition:  
there is a Seifert hypersurface $V'_i$ for $K'_i$ ($i=1,2$) 
such that $V'_1\cap V'_2$ is one 2-sphere $S^2_0$. 
}
\vskip3mm

\noindent{\bf Proof of Lemma 8.1.2. }   
Take $V_1$ and $V_2$ in 
Lemma 8.1.1.   
If $\nu=0$, then 
Theorem 8.1   
holds. 
If $\nu=1$, then 
Lemma 8.1.2   
holds. 
Suppose $\nu\geq2$. 

We can suppose that $S^2_1$ and $S^2_2$ satisfy the following: 
There is a point $p_1\in S^2_1$, a point $p_2\in S^2_2$, and 
a path $l\subset V_1$ such that 
(1) $\partial l=p_1\amalg p_2$ 
(2) $l\cap (S^2_1\amalg...\amalg S^2_\nu)=p_1\amalg p_2$  
(3) $l\cap K_1=\phi.$

Take a 4-dimensional 1-handle $h^1\subset S^4$ whose core is $l$ 
such that $h^1$ is attached to $V_2$ along $p_1\amalg p_2$. 
Then $h^1\cap V_1$ is a 3-dimensional 1-handle 
which is attached to $S^2_1\amalg S^2_2$ along $p_1\amalg p_2$. 
We carry out surgery on $V_2$ by using $h^1$. 
The new manifold is called $V^\natural_2$.
Then $V^\natural_2$ is a connected Seifert hypersurface for $K_2$. 
When we carry out the surgery on $V_2$, 
we carry out surgery on $S^2_1\amalg S^2_2$ 
by using the 3-dimensional 1-handle $h^1\cap V_1$. 
Then the result is a 2-sphere. 
Then $V_1\cap V^\natural_2$ is $(\nu-1)$ 2-spheres. 
By the induction on $\nu$, 
Lemma 8.1.2   
holds.

\vskip3mm
\noindent
{\bf Lemma 8.1.3.}   
 {\it
Let $L=(K_1,K_2)$ be the 2-link in 
Theorem 8.1.   
Then there is a 2-link $L''=(K''_1,K''_2)$ which is equivalent to $L$ 
satisfying the following condition:  
there is a connected Seifert hypersurface $V''_i$ for $K''_i$ ($i=1,2$) 
such that $V''_1\cap V''_2$ is one 2-disc $D^2_0$. 
}
\vskip3mm

\noindent{\bf Proof of Lemma 8.1.3. }   
Take $V'_1$ and $V'_2$ in 
Lemma 8.1.2.   
Take a point $p\subset K'_1=\partial V'_1$. 
Take a point $q\subset S^2_0=V'_1\cap V'_2$. 
Take a path $l\subset V'_1$ such that 
(1) $\partial l=p\amalg q$. 
(2) $l\cap S^2_0=q$  
(3) $l\cap K'_1=p.$

Let $N$ be a tubular neighborhood of $l$ in $V'_1$, 
Then $N$ is a 3-ball. 
Note that $N\cap K_1$ is a 2-disc, which is a tubular neighborhood of 
$p$ in $K_1$. 
Note that $q\subset $Int $N$. 
Note that Int ($N\cap S^2_0$) is in Int $N$.   
Then the following hold. 
(1) $\overline{V'_1-N}$ $\cap V'_2$ is a 2-disc. 
(2) $\partial(\overline{V'_1-N})$ is equivalent to $K_1$. 
(3) $(\partial(\overline{V'_1-N}), K_2)$ is equivalent to $L'$ 
and hence to $L$. 
$(\partial(\overline{V'_1-N}), K_2)$ is called $L''=(K''_1,K''_2)$. 
This completes the proof of 
Lemma 8.1.3.    

\noindent{\bf Proof of Theorem 8.1. }   
Take $V''_1$ and $V''_2$ in 
Lemma 8.1.3.   
We can suppose that $\partial D^2_0\subset K''_1$ and that 
$D^2\subset$ Int$ V''_2$. 
Take a 3-ball $P\subset V''_2$ such that 
$P\cap K''_2$ is a 2-disc and that 
$D^2\subset$ Int $P$. 
Then $V''_1\cap \partial(\overline{V''_2-P})=\phi$. 
Let $L^!$ be a 2-link ($ K''_1, \partial(\overline{V''_2-P}) $). 
Then  $L^!$ is a boundary link. 
Furthermore $L^!$ is obtained from $L''$ by 
an operation that we fix $K''_1$ and 
that we move $K''_2$ to $\partial(\overline{V''_2-P})$ 
so that we fix $K''_2\cap\partial(\overline{V''_2-P})$.  
This operation on $L^!$ is essentially same as a ribbon move. 
This completes the proof of 
Theorem 8.1.   3


\vskip3mm
\noindent
{\bf Theorem 8.2}   
{\it
Let $L=(K_1,K_2)$ be an SHB 2-link. Then $\mu(L)=\mu(K_1)+\mu(K_2)$. 
}
\vskip3mm

\noindent{\bf Proof.}  
By Theorem 8.1,   
$L$ is ribbon-move equivalent to a boundary 2-link 
$\bar{L}=(\bar{K_1}, \bar{K_2})$.  
Let 
$\bar{V_i}$ be a Seifert hypersurface for  
$\bar{K_i}$ such that $\bar{V_1}\cap\bar{V_2}=\phi$. 
Then $\mu(\bar{L})$=$\mu(\bar{V_1}\cup h^3)+\mu(\bar{V_2}\cup h^3)$, 
where $h^3$ is a 3-dimensional 3-handle which is attached to 
$\bar{V_i}$ along 
the 2-sphere 
$\partial\bar{V_i}$.
Hence $\mu(\bar{L})$=$\mu(\bar{K_1})+\mu(\bar{K_2})$. 
By Theorem 2.2,   
$\mu(L)=\mu(\bar{L})$ and $\mu(K_i)=\mu(\bar{K_i})$. 
Hence $\mu(L)=\mu(K_1)+\mu(K_2)$. 
\

\vskip3mm \noindent{\bf Problem 8.3.}   

(1)  
Let $L=(K_1,K_2)$ be a 2-link. Then does $\mu(L)=\mu(K_1)+\mu(K_2)$ hold? 

(2)   
Is there an $n$-link which is not an SHB link ($n\geq2$)?

\vskip3mm
  \noindent{\large{\bf\S9.    
Discussions
 }}

We would point out the following facts 
by analogy of the discussions of finite type invariants of 1-knots 
(e.g. \cite{BL}) 
although they are very easy observations.

By using Theorem 2.2 we have:  
The $\mu$-invariant of 2-links is an order zero finite type invariant 
if we define `order of invariants' by using ribbon-moves
( e.g. as follows ),   
and 
there is a 2-knot whose $\mu$-invariant is not zero.

We define order, for example, as follows. 
Let $I_n$ be the set of immersed $m$ 2-spheres 
with the conditions: 
(1) The set of singular points consists of double points.  
(2) Each component of the set of singular points is as in Figure 9.3. 
(3) The components of the set of singular points are $n$. 
Then $I_0$ is the set of $m$-component 2-links. 
Let $v_i(\quad)\in G$ 
be an invariant of elements of $I_i$, 
where $G$ is a group. 
Let $X_0$ be an element of $I_{i+1}$. 
Let $X_+$ and $X_-$ be elements of $I_{i}$. 
Suppose that 
$X_0$, $X_+$ and $X_-$ coincide in $S^4-B^4$.  
Suppose that 
$X_0\cap B$ is drawn as in Figure 9.3, 
$X_+\cap B$ is drawn as in Figure 9.1, and 
$X_-\cap B$ is drawn as in Figure 9.2. 
In Figure 9.1, 9.2, 9.3, 
we do not assume the orientation of $X_*\cap B$ and that of $B$.  
If we have 
$\{v_{i+1}(X)\}^2=\{v_i(X_+)-v_i(X_-)\}^2$ 
and $v_i$ is zero for $i>p$, 
then we call $v_*(\quad)$ is an order $p$ invariant of 2-links.

We define a link-type invariant $v(\quad)$ of $(S^2, T^2)$-links. 
( See 
\cite{S}  
for detail. 
See $(S^2, T^2)$-links 
for 
\cite{O1}.)
We call it the alinking number of $(S^2, T^2)$-links. 
Let $L=(L_S, L_T)$ be a $(S^2, T^2)$-link. 
Let $\iota$ be the map 
$H^1(S^4-L_S;{\bf Z})\rightarrow H^1(L_T;{\bf Z})$ induced by the inclusion. 


Define 
\[
v(L)=
\left\{ 
\begin{array}{ll}
n & \mbox{if $H^1(L_T;{\bf Z})$/Im$\iota\cong{\bf Z}\oplus({\bf Z}/(n\cdot{\bf Z}))$ $(n\geq2, n\in {\bf N})$}\\ 
1& \mbox{if $H^1(L_T;{\bf Z})$/Im$\iota\cong{\bf Z}$}\\
0& \mbox{if $H^1(L_T;{\bf Z})$/Im$\iota\cong{\bf Z}\oplus{\bf Z}$.}
\end{array} 
\right.
\]

Then the mod 2 alinking number of $(S^2, T^2)$-links is 
an order one finite type invariant 
if we define `order of invariants' by using ribbon-moves
(e.g. as above),  
and 
there is an $(S^2, T^2)$-link whose mod 2 alinking number is not zero.  
(The proof is similar to the proof that 
the linking number of 2-component 1-links  
is an order one finite type invariant. See
\cite{BL}.)  

\noindent{\bf Note.}
\cite{B}  
 and 
\cite{CCFM} 
 etc. try to make a high-dimensional version of 
works on 1-links 
by Jones, Witten, Kontsevich, Vassiliev, etc.  
(in 
\cite{J}    
\cite{Ko}    
\cite{V}    
\cite{W}    
etc.
)

\footnotesize{
 }

\noindent{\large{\bf Appendix }}

{\normalsize 
A { \it ribbon 2-link } is a 2-link 
$L=(K_1,...,K_m)$ with the following properties. 
There is a self-transverse immersion 
$f:D^3_1\amalg...\amalg D^3_m\to S^4$ such that: 
(a)$ f(\partial D^3_i)$ coincides with $K_i$. 
(b)The singular point set $X$ consists of double points.  
(c)For each connected component $X_i$ of $X$, 
$f^{-1}(X_i)$ is diffeomorphic to the two 2-discs. 
(d)Put $\partial\{f^{-1}(X_i)\}$=$P\amalg Q$. 
One of  $P\amalg Q$ is included in the boundary of $D^3_i$ and  
another of $P\amalg Q$ is included in the interior of $D^3_j$ 
for integers $i, j$.
( We do not assume $i\neq j$ nor $i=j$. ) 
}

\np 

\pagestyle{empty}

\unitlength 0.1in
\begin{picture}(56.10,47.00)(8.50,-47.60)
%
\special{pn 8}%
\special{ar 3510 320 560 250  0.0000000 6.2831853}%
%
\special{pn 20}%
\special{pa 3660 3950}%
\special{pa 3660 2730}%
\special{fp}%
%
\special{pn 20}%
\special{pa 3340 3970}%
\special{pa 3340 2740}%
\special{fp}%
%
\special{pn 20}%
\special{pa 3350 4000}%
\special{pa 3360 3970}%
\special{pa 3384 3949}%
\special{pa 3413 3935}%
\special{pa 3444 3926}%
\special{pa 3475 3921}%
\special{pa 3507 3920}%
\special{pa 3539 3923}%
\special{pa 3570 3930}%
\special{pa 3601 3940}%
\special{pa 3627 3958}%
\special{pa 3647 3983}%
\special{pa 3648 4014}%
\special{pa 3630 4040}%
\special{pa 3603 4058}%
\special{pa 3573 4070}%
\special{pa 3542 4077}%
\special{pa 3510 4080}%
\special{pa 3478 4079}%
\special{pa 3447 4075}%
\special{pa 3416 4066}%
\special{pa 3387 4053}%
\special{pa 3362 4032}%
\special{pa 3350 4003}%
\special{pa 3350 4000}%
\special{sp}%
%
\special{pn 20}%
\special{pa 3340 360}%
\special{pa 3340 1580}%
\special{fp}%
%
\special{pn 20}%
\special{pa 3660 340}%
\special{pa 3660 1570}%
\special{fp}%
%
\special{pn 20}%
\special{ar 3500 310 150 80  0.0000000 6.2831853}%
%
\special{pn 8}%
\special{ar 3510 4080 560 250  0.0000000 6.2831853}%
%
\special{pn 8}%
\special{pa 4080 330}%
\special{pa 4080 4080}%
\special{fp}%
%
\special{pn 8}%
\special{pa 2950 340}%
\special{pa 2950 4080}%
\special{fp}%
%
\special{pn 8}%
\special{ar 1410 330 560 250  0.0000000 6.2831853}%
%
\special{pn 8}%
\special{ar 1410 4090 560 250  0.0000000 6.2831853}%
%
\special{pn 8}%
\special{pa 1980 340}%
\special{pa 1980 4090}%
\special{fp}%
%
\special{pn 8}%
\special{pa 850 350}%
\special{pa 850 4090}%
\special{fp}%
%
\special{pn 8}%
\special{ar 5890 310 560 250  0.0000000 6.2831853}%
%
\special{pn 8}%
\special{ar 5890 4070 560 250  0.0000000 6.2831853}%
%
\special{pn 8}%
\special{pa 6460 320}%
\special{pa 6460 4070}%
\special{fp}%
%
\special{pn 8}%
\special{pa 5330 330}%
\special{pa 5330 4070}%
\special{fp}%
%
\special{pn 20}%
\special{ar 3500 1590 150 80  0.0000000 6.2831853}%
%
\special{pn 20}%
\special{ar 3500 2710 150 80  0.0000000 6.2831853}%
%
\special{pn 20}%
\special{ar 3520 2070 560 250  0.0000000 6.2831853}%
%
\special{pn 8}%
\special{pa 3660 1600}%
\special{pa 5730 1600}%
\special{dt 0.045}%
\special{pa 5730 1600}%
\special{pa 5729 1600}%
\special{dt 0.045}%
%
\special{pn 8}%
\special{pa 3680 2720}%
\special{pa 5750 2720}%
\special{dt 0.045}%
\special{pa 5750 2720}%
\special{pa 5749 2720}%
\special{dt 0.045}%
%
\special{pn 20}%
\special{ar 5850 1600 150 80  0.0000000 6.2831853}%
%
\special{pn 20}%
\special{ar 5850 2730 150 80  0.0000000 6.2831853}%
%
\special{pn 20}%
\special{pa 5690 1610}%
\special{pa 5690 2700}%
\special{fp}%
%
\special{pn 20}%
\special{pa 6020 1640}%
\special{pa 6020 2730}%
\special{fp}%
\put(12.0000,-46.0000){\makebox(0,0)[lb]{t=-0.5}}%
\put(32.0000,-46.0000){\makebox(0,0)[lb]{t=0}}%
\put(56.0000,-46.0000){\makebox(0,0)[lb]{t=0.5}}%
\put(30.3000,-49.3000){\makebox(0,0)[lb]{Figure 1.1}}%
%
\special{pn 8}%
\special{ar 3510 320 560 250  0.0000000 6.2831853}%
%
\special{pn 8}%
\special{pa 2950 340}%
\special{pa 2950 4080}%
\special{fp}%
%
\special{pn 8}%
\special{pa 4080 330}%
\special{pa 4080 4080}%
\special{fp}%
%
\special{pn 8}%
\special{ar 3510 4080 560 250  0.0000000 6.2831853}%
%
\special{pn 20}%
\special{ar 3500 310 150 80  0.0000000 6.2831853}%
%
\special{pn 20}%
\special{pa 3660 340}%
\special{pa 3660 1570}%
\special{fp}%
%
\special{pn 20}%
\special{pa 3340 360}%
\special{pa 3340 1580}%
\special{fp}%
%
\special{pn 20}%
\special{pa 3350 4000}%
\special{pa 3360 3970}%
\special{pa 3384 3949}%
\special{pa 3413 3935}%
\special{pa 3444 3926}%
\special{pa 3475 3921}%
\special{pa 3507 3920}%
\special{pa 3539 3923}%
\special{pa 3570 3930}%
\special{pa 3601 3940}%
\special{pa 3627 3958}%
\special{pa 3647 3983}%
\special{pa 3648 4014}%
\special{pa 3630 4040}%
\special{pa 3603 4058}%
\special{pa 3573 4070}%
\special{pa 3542 4077}%
\special{pa 3510 4080}%
\special{pa 3478 4079}%
\special{pa 3447 4075}%
\special{pa 3416 4066}%
\special{pa 3387 4053}%
\special{pa 3362 4032}%
\special{pa 3350 4003}%
\special{pa 3350 4000}%
\special{sp}%
%
\special{pn 20}%
\special{pa 3340 3970}%
\special{pa 3340 2740}%
\special{fp}%
%
\special{pn 20}%
\special{pa 3660 3950}%
\special{pa 3660 2730}%
\special{fp}%
%
\special{pn 8}%
\special{ar 1410 330 560 250  0.0000000 6.2831853}%
%
\special{pn 8}%
\special{ar 1410 4090 560 250  0.0000000 6.2831853}%
%
\special{pn 8}%
\special{pa 1980 340}%
\special{pa 1980 4090}%
\special{fp}%
%
\special{pn 8}%
\special{pa 850 350}%
\special{pa 850 4090}%
\special{fp}%
%
\special{pn 8}%
\special{ar 5890 310 560 250  0.0000000 6.2831853}%
%
\special{pn 8}%
\special{ar 5890 4070 560 250  0.0000000 6.2831853}%
%
\special{pn 8}%
\special{pa 6460 320}%
\special{pa 6460 4070}%
\special{fp}%
%
\special{pn 8}%
\special{pa 5330 330}%
\special{pa 5330 4070}%
\special{fp}%
\end{picture}%
 
\np 
\unitlength 0.1in
\begin{picture}(56.10,47.10)(8.50,-47.70)
%
\special{pn 8}%
\special{ar 3510 320 560 250  0.0000000 6.2831853}%
%
\special{pn 20}%
\special{pa 3660 3950}%
\special{pa 3660 2730}%
\special{fp}%
%
\special{pn 20}%
\special{pa 3340 3970}%
\special{pa 3340 2740}%
\special{fp}%
%
\special{pn 20}%
\special{pa 3350 4000}%
\special{pa 3360 3970}%
\special{pa 3384 3949}%
\special{pa 3413 3935}%
\special{pa 3444 3926}%
\special{pa 3475 3921}%
\special{pa 3507 3920}%
\special{pa 3539 3923}%
\special{pa 3570 3930}%
\special{pa 3601 3940}%
\special{pa 3627 3958}%
\special{pa 3647 3983}%
\special{pa 3648 4014}%
\special{pa 3630 4040}%
\special{pa 3603 4058}%
\special{pa 3573 4070}%
\special{pa 3542 4077}%
\special{pa 3510 4080}%
\special{pa 3478 4079}%
\special{pa 3447 4075}%
\special{pa 3416 4066}%
\special{pa 3387 4053}%
\special{pa 3362 4032}%
\special{pa 3350 4003}%
\special{pa 3350 4000}%
\special{sp}%
%
\special{pn 20}%
\special{pa 3340 360}%
\special{pa 3340 1580}%
\special{fp}%
%
\special{pn 20}%
\special{pa 3660 340}%
\special{pa 3660 1570}%
\special{fp}%
%
\special{pn 20}%
\special{ar 3500 310 150 80  0.0000000 6.2831853}%
%
\special{pn 8}%
\special{ar 3510 4080 560 250  0.0000000 6.2831853}%
%
\special{pn 8}%
\special{pa 4080 330}%
\special{pa 4080 4080}%
\special{fp}%
%
\special{pn 8}%
\special{pa 2950 340}%
\special{pa 2950 4080}%
\special{fp}%
%
\special{pn 8}%
\special{ar 1410 330 560 250  0.0000000 6.2831853}%
%
\special{pn 8}%
\special{ar 1410 4090 560 250  0.0000000 6.2831853}%
%
\special{pn 8}%
\special{pa 1980 340}%
\special{pa 1980 4090}%
\special{fp}%
%
\special{pn 8}%
\special{pa 850 350}%
\special{pa 850 4090}%
\special{fp}%
%
\special{pn 8}%
\special{ar 5890 310 560 250  0.0000000 6.2831853}%
%
\special{pn 8}%
\special{ar 5890 4070 560 250  0.0000000 6.2831853}%
%
\special{pn 8}%
\special{pa 6460 320}%
\special{pa 6460 4070}%
\special{fp}%
%
\special{pn 8}%
\special{pa 5330 330}%
\special{pa 5330 4070}%
\special{fp}%
%
\special{pn 20}%
\special{ar 3500 1590 150 80  0.0000000 6.2831853}%
%
\special{pn 20}%
\special{ar 3500 2710 150 80  0.0000000 6.2831853}%
%
\special{pn 20}%
\special{ar 3520 2070 560 250  0.0000000 6.2831853}%
\put(12.0000,-46.0000){\makebox(0,0)[lb]{t=-0.5}}%
\put(32.0000,-46.0000){\makebox(0,0)[lb]{t=0}}%
\put(56.0000,-46.0000){\makebox(0,0)[lb]{t=0.5}}%
%
\special{pn 8}%
\special{ar 3510 320 560 250  0.0000000 6.2831853}%
%
\special{pn 8}%
\special{pa 2950 340}%
\special{pa 2950 4080}%
\special{fp}%
%
\special{pn 8}%
\special{pa 4080 330}%
\special{pa 4080 4080}%
\special{fp}%
%
\special{pn 8}%
\special{ar 3510 4080 560 250  0.0000000 6.2831853}%
%
\special{pn 20}%
\special{ar 3500 310 150 80  0.0000000 6.2831853}%
%
\special{pn 20}%
\special{pa 3660 340}%
\special{pa 3660 1570}%
\special{fp}%
%
\special{pn 20}%
\special{pa 3340 360}%
\special{pa 3340 1580}%
\special{fp}%
%
\special{pn 20}%
\special{pa 3350 4000}%
\special{pa 3360 3970}%
\special{pa 3384 3949}%
\special{pa 3413 3935}%
\special{pa 3444 3926}%
\special{pa 3475 3921}%
\special{pa 3507 3920}%
\special{pa 3539 3923}%
\special{pa 3570 3930}%
\special{pa 3601 3940}%
\special{pa 3627 3958}%
\special{pa 3647 3983}%
\special{pa 3648 4014}%
\special{pa 3630 4040}%
\special{pa 3603 4058}%
\special{pa 3573 4070}%
\special{pa 3542 4077}%
\special{pa 3510 4080}%
\special{pa 3478 4079}%
\special{pa 3447 4075}%
\special{pa 3416 4066}%
\special{pa 3387 4053}%
\special{pa 3362 4032}%
\special{pa 3350 4003}%
\special{pa 3350 4000}%
\special{sp}%
%
\special{pn 20}%
\special{pa 3340 3970}%
\special{pa 3340 2740}%
\special{fp}%
%
\special{pn 20}%
\special{pa 3660 3950}%
\special{pa 3660 2730}%
\special{fp}%
%
\special{pn 8}%
\special{ar 1410 330 560 250  0.0000000 6.2831853}%
%
\special{pn 8}%
\special{ar 1410 4090 560 250  0.0000000 6.2831853}%
%
\special{pn 8}%
\special{pa 1980 340}%
\special{pa 1980 4090}%
\special{fp}%
%
\special{pn 8}%
\special{pa 850 350}%
\special{pa 850 4090}%
\special{fp}%
%
\special{pn 8}%
\special{ar 5890 310 560 250  0.0000000 6.2831853}%
%
\special{pn 8}%
\special{ar 5890 4070 560 250  0.0000000 6.2831853}%
%
\special{pn 8}%
\special{pa 6460 320}%
\special{pa 6460 4070}%
\special{fp}%
%
\special{pn 8}%
\special{pa 5330 330}%
\special{pa 5330 4070}%
\special{fp}%
%
\special{pn 8}%
\special{pa 3320 1590}%
\special{pa 1590 1590}%
\special{dt 0.045}%
\special{pa 1590 1590}%
\special{pa 1591 1590}%
\special{dt 0.045}%
%
\special{pn 8}%
\special{pa 3330 2740}%
\special{pa 1620 2740}%
\special{dt 0.045}%
\special{pa 1620 2740}%
\special{pa 1621 2740}%
\special{dt 0.045}%
%
\special{pn 20}%
\special{ar 1420 1590 150 80  0.0000000 6.2831853}%
%
\special{pn 20}%
\special{ar 1420 2720 150 80  0.0000000 6.2831853}%
%
\special{pn 20}%
\special{pa 1260 1600}%
\special{pa 1260 2690}%
\special{fp}%
%
\special{pn 20}%
\special{pa 1590 1630}%
\special{pa 1590 2720}%
\special{fp}%
\put(30.0000,-49.4000){\makebox(0,0)[lb]{Figure 1.2}}%
\end{picture}%
 
\np 
\input 4.1.tex

\np 
\input 4.2.tex

\np\rotatebox[origin=c]{90}{
\unitlength 0.1in
\begin{picture}(64.20,45.60)(0.90,-46.00)
%
\special{pn 8}%
\special{ar 2430 290 560 250  0.0000000 6.2831853}%
\put(62.3000,-47.7000){\makebox(0,0)[lb]{Figure 4.3}}%
%
\special{pn 8}%
\special{pa 4080 1120}%
\special{pa 4080 3250}%
\special{fp}%
%
\special{pn 8}%
\special{pa 4370 1110}%
\special{pa 4370 3230}%
\special{fp}%
%
\special{pn 8}%
\special{pa 2570 3290}%
\special{pa 4060 3290}%
\special{dt 0.045}%
\special{pa 4060 3290}%
\special{pa 4059 3290}%
\special{dt 0.045}%
%
\special{pn 8}%
\special{pa 2590 1050}%
\special{pa 4070 1050}%
\special{dt 0.045}%
\special{pa 4070 1050}%
\special{pa 4069 1050}%
\special{dt 0.045}%
%
\special{pn 8}%
\special{pa 2270 3950}%
\special{pa 2270 3300}%
\special{fp}%
%
\special{pn 8}%
\special{pa 2570 3950}%
\special{pa 2570 3300}%
\special{fp}%
%
\special{pn 8}%
\special{pa 2270 330}%
\special{pa 2270 1000}%
\special{fp}%
%
\special{pn 8}%
\special{pa 2570 310}%
\special{pa 2570 980}%
\special{fp}%
\put(56.4000,-45.7000){\makebox(0,0)[lb]{t=0.7}}%
%
\special{pn 8}%
\special{ar 5940 4050 560 250  0.0000000 6.2831853}%
%
\special{pn 8}%
\special{pa 6510 300}%
\special{pa 6510 4050}%
\special{fp}%
%
\special{pn 8}%
\special{pa 5380 310}%
\special{pa 5380 4050}%
\special{fp}%
%
\special{pn 8}%
\special{ar 5940 290 560 250  0.0000000 6.2831853}%
%
\special{pn 8}%
\special{ar 5940 4050 560 250  0.0000000 6.2831853}%
%
\special{pn 8}%
\special{pa 6510 300}%
\special{pa 6510 4050}%
\special{fp}%
%
\special{pn 8}%
\special{pa 5380 310}%
\special{pa 5380 4050}%
\special{fp}%
%
\special{pn 8}%
\special{ar 5940 290 560 250  0.0000000 6.2831853}%
%
\special{pn 8}%
\special{pa 90 320}%
\special{pa 90 4060}%
\special{fp}%
%
\special{pn 8}%
\special{pa 1220 310}%
\special{pa 1220 4060}%
\special{fp}%
%
\special{pn 8}%
\special{ar 650 4060 560 250  0.0000000 6.2831853}%
%
\special{pn 8}%
\special{ar 650 300 560 250  0.0000000 6.2831853}%
\put(4.4000,-45.7000){\makebox(0,0)[lb]{t=-0.5}}%
%
\special{pn 8}%
\special{ar 4270 4070 560 250  0.0000000 6.2831853}%
%
\special{pn 8}%
\special{pa 4840 320}%
\special{pa 4840 4070}%
\special{fp}%
%
\special{pn 8}%
\special{pa 3710 330}%
\special{pa 3710 4070}%
\special{fp}%
%
\special{pn 8}%
\special{ar 4230 1050 150 80  0.0000000 6.2831853}%
%
\special{pn 8}%
\special{ar 4230 3290 150 80  0.0000000 6.2831853}%
\put(39.8000,-46.0000){\makebox(0,0)[lb]{t=0.5}}%
%
\special{pn 8}%
\special{ar 4270 310 560 250  0.0000000 6.2831853}%
%
\special{pn 8}%
\special{ar 4270 4070 560 250  0.0000000 6.2831853}%
%
\special{pn 8}%
\special{pa 4840 320}%
\special{pa 4840 4070}%
\special{fp}%
%
\special{pn 8}%
\special{pa 3710 330}%
\special{pa 3710 4070}%
\special{fp}%
%
\special{pn 8}%
\special{ar 4270 310 560 250  0.0000000 6.2831853}%
%
\special{pn 8}%
\special{pa 90 320}%
\special{pa 90 4060}%
\special{fp}%
%
\special{pn 8}%
\special{pa 1220 310}%
\special{pa 1220 4060}%
\special{fp}%
%
\special{pn 8}%
\special{ar 650 4060 560 250  0.0000000 6.2831853}%
%
\special{pn 8}%
\special{ar 650 300 560 250  0.0000000 6.2831853}%
%
\special{pn 8}%
\special{pa 3000 300}%
\special{pa 3000 4050}%
\special{fp}%
%
\special{pn 8}%
\special{pa 1870 310}%
\special{pa 1870 4050}%
\special{fp}%
%
\special{pn 8}%
\special{ar 2420 1030 150 80  0.0000000 6.2831853}%
%
\special{pn 8}%
\special{ar 2410 3280 150 80  0.0000000 6.2831853}%
\put(21.2000,-45.7000){\makebox(0,0)[lb]{t=0}}%
%
\special{pn 8}%
\special{ar 2430 290 560 250  0.0000000 6.2831853}%
%
\special{pn 8}%
\special{pa 1870 310}%
\special{pa 1870 4050}%
\special{fp}%
%
\special{pn 8}%
\special{pa 3000 300}%
\special{pa 3000 4050}%
\special{fp}%
%
\special{pn 8}%
\special{ar 2430 4050 560 250  0.0000000 6.2831853}%
%
\special{pn 8}%
\special{ar 2420 280 150 80  0.0000000 6.2831853}%
%
\special{pn 8}%
\special{pa 2270 3970}%
\special{pa 2280 3940}%
\special{pa 2304 3919}%
\special{pa 2333 3905}%
\special{pa 2364 3896}%
\special{pa 2395 3891}%
\special{pa 2427 3890}%
\special{pa 2459 3893}%
\special{pa 2490 3900}%
\special{pa 2521 3910}%
\special{pa 2547 3928}%
\special{pa 2567 3953}%
\special{pa 2568 3984}%
\special{pa 2550 4010}%
\special{pa 2523 4028}%
\special{pa 2493 4040}%
\special{pa 2462 4047}%
\special{pa 2430 4050}%
\special{pa 2398 4049}%
\special{pa 2367 4045}%
\special{pa 2336 4036}%
\special{pa 2307 4023}%
\special{pa 2282 4002}%
\special{pa 2270 3973}%
\special{pa 2270 3970}%
\special{sp}%
%
\special{pn 8}%
\special{ar 2440 1430 560 250  0.0000000 6.2831853}%
\end{picture}%
}

 
\np\rotatebox[origin=c]{90}{\input {4.4.tex}}


\np 
\unitlength 0.1in
\begin{picture}(46.40,46.60)(17.20,-47.80)
%
\special{pn 8}%
\special{ar 2280 370 560 250  0.0000000 6.2831853}%
\put(40.1000,-47.5000){\makebox(0,0)[lb]{t=0}}%
%
\special{pn 8}%
\special{pa 6360 380}%
\special{pa 6360 4130}%
\special{fp}%
%
\special{pn 8}%
\special{pa 5230 390}%
\special{pa 5230 4130}%
\special{fp}%
%
\special{pn 8}%
\special{ar 5790 370 560 250  0.0000000 6.2831853}%
%
\special{pn 8}%
\special{ar 5790 4130 560 250  0.0000000 6.2831853}%
%
\special{pn 8}%
\special{pa 6360 380}%
\special{pa 6360 4130}%
\special{fp}%
%
\special{pn 8}%
\special{pa 5230 390}%
\special{pa 5230 4130}%
\special{fp}%
%
\special{pn 8}%
\special{ar 5790 370 560 250  0.0000000 6.2831853}%
\put(20.1000,-47.4000){\makebox(0,0)[lb]{t=-0.5}}%
%
\special{pn 8}%
\special{ar 4120 4150 560 250  0.0000000 6.2831853}%
%
\special{pn 8}%
\special{pa 4690 400}%
\special{pa 4690 4150}%
\special{fp}%
%
\special{pn 8}%
\special{pa 3560 410}%
\special{pa 3560 4150}%
\special{fp}%
%
\special{pn 13}%
\special{ar 4080 410 150 80  0.0000000 6.2831853}%
%
\special{pn 13}%
\special{ar 4080 4120 150 80  0.0000000 6.2831853}%
\put(55.6000,-47.6000){\makebox(0,0)[lb]{t=0.5}}%
%
\special{pn 8}%
\special{ar 4120 390 560 250  0.0000000 6.2831853}%
%
\special{pn 8}%
\special{ar 4120 4150 560 250  0.0000000 6.2831853}%
%
\special{pn 8}%
\special{pa 4690 400}%
\special{pa 4690 4150}%
\special{fp}%
%
\special{pn 8}%
\special{pa 3560 410}%
\special{pa 3560 4150}%
\special{fp}%
%
\special{pn 8}%
\special{ar 4120 390 560 250  0.0000000 6.2831853}%
%
\special{pn 8}%
\special{ar 2280 4130 560 250  0.0000000 6.2831853}%
%
\special{pn 8}%
\special{pa 2850 380}%
\special{pa 2850 4130}%
\special{fp}%
%
\special{pn 8}%
\special{pa 1720 390}%
\special{pa 1720 4130}%
\special{fp}%
%
\special{pn 8}%
\special{ar 2280 370 560 250  0.0000000 6.2831853}%
%
\special{pn 8}%
\special{pa 1720 390}%
\special{pa 1720 4130}%
\special{fp}%
%
\special{pn 8}%
\special{pa 2850 380}%
\special{pa 2850 4130}%
\special{fp}%
%
\special{pn 8}%
\special{ar 2280 4130 560 250  0.0000000 6.2831853}%
%
\special{pn 13}%
\special{pa 4230 2040}%
\special{pa 4230 420}%
\special{fp}%
\put(33.0000,-49.5000){\makebox(0,0)[lb]{Figure 4.5}}%
%
\special{pn 13}%
\special{pa 3920 4090}%
\special{pa 3920 430}%
\special{fp}%
%
\special{pn 13}%
\special{pa 4230 440}%
\special{pa 4230 4120}%
\special{fp}%
\end{picture}%

\np 
\input 4.6.tex

\np
\rotatebox[origin=c]{90}{
\unitlength 0.1in
\begin{picture}(73.38,24.30)(3.70,-34.70)
%
\special{pn 8}%
\special{pa 7174 1956}%
\special{pa 6254 1956}%
\special{pa 6267 1040}%
\special{pa 7708 1040}%
\special{pa 7708 1956}%
\special{pa 7700 2484}%
\special{pa 7121 2484}%
\special{pa 7121 2668}%
\special{pa 7121 2668}%
\special{fp}%
%
\special{pn 4}%
\special{sh 0}%
\special{pa 6885 1008}%
\special{pa 6944 1008}%
\special{pa 6944 1080}%
\special{pa 6885 1080}%
\special{pa 6885 1008}%
\special{ip}%
%
\special{pn 4}%
\special{sh 0}%
\special{pa 7030 1015}%
\special{pa 7089 1015}%
\special{pa 7089 1086}%
\special{pa 7030 1086}%
\special{pa 7030 1015}%
\special{ip}%
%
\special{pn 4}%
\special{sh 0}%
\special{pa 6627 1924}%
\special{pa 6689 1924}%
\special{pa 6689 1996}%
\special{pa 6627 1996}%
\special{pa 6627 1924}%
\special{ip}%
%
\special{pn 4}%
\special{sh 0}%
\special{pa 6484 1931}%
\special{pa 6543 1931}%
\special{pa 6543 2003}%
\special{pa 6484 2003}%
\special{pa 6484 1931}%
\special{ip}%
%
\special{pn 8}%
\special{pa 6910 639}%
\special{pa 6918 1555}%
\special{pa 6516 1555}%
\special{pa 6516 2346}%
\special{pa 6971 2346}%
\special{pa 6971 2668}%
\special{pa 6971 2668}%
\special{pa 6971 2668}%
\special{fp}%
%
\special{pn 8}%
\special{pa 7049 639}%
\special{pa 7049 1694}%
\special{pa 6648 1694}%
\special{pa 6648 2220}%
\special{pa 6965 2220}%
\special{pa 6971 2088}%
\special{pa 7182 2088}%
\special{pa 7174 1956}%
\special{pa 7174 1956}%
\special{pa 7174 1956}%
\special{fp}%
%
\special{pn 8}%
\special{pa 5997 639}%
\special{pa 7990 639}%
\special{pa 7990 2674}%
\special{pa 5997 2674}%
\special{pa 5997 639}%
\special{dt 0.045}%
%
\special{pn 8}%
\special{pa 4804 2503}%
\special{pa 4804 2674}%
\special{fp}%
%
\special{pn 8}%
\special{pa 4962 1963}%
\special{pa 4041 1963}%
\special{pa 4055 1046}%
\special{pa 5494 1046}%
\special{pa 5494 1963}%
\special{pa 5489 2489}%
\special{pa 4910 2489}%
\special{pa 4910 2674}%
\special{pa 4910 2674}%
\special{fp}%
%
\special{pn 8}%
\special{pa 3785 646}%
\special{pa 5778 646}%
\special{pa 5778 2681}%
\special{pa 3785 2681}%
\special{pa 3785 646}%
\special{dt 0.045}%
%
\special{pn 4}%
\special{sh 0}%
\special{pa 4673 1013}%
\special{pa 4869 1013}%
\special{pa 4869 1093}%
\special{pa 4673 1093}%
\special{pa 4673 1013}%
\special{ip}%
%
\special{pn 4}%
\special{sh 0}%
\special{pa 4266 1938}%
\special{pa 4463 1938}%
\special{pa 4463 2015}%
\special{pa 4266 2015}%
\special{pa 4266 1938}%
\special{ip}%
%
\special{pn 8}%
\special{pa 4837 646}%
\special{pa 4837 1698}%
\special{pa 4437 1698}%
\special{pa 4437 2226}%
\special{pa 4751 2226}%
\special{pa 4758 2095}%
\special{pa 4969 2095}%
\special{pa 4962 1963}%
\special{pa 4962 1963}%
\special{pa 4962 1963}%
\special{fp}%
%
\special{pn 8}%
\special{pa 4699 646}%
\special{pa 4707 1561}%
\special{pa 4305 1561}%
\special{pa 4305 2353}%
\special{pa 4758 2353}%
\special{pa 4758 2674}%
\special{pa 4758 2674}%
\special{pa 4758 2674}%
\special{fp}%
%
\special{pn 8}%
\special{pa 4785 646}%
\special{pa 4785 1694}%
\special{fp}%
%
\special{pn 8}%
\special{pa 4739 646}%
\special{pa 4739 1687}%
\special{fp}%
%
\special{pn 8}%
\special{pa 4673 1561}%
\special{pa 4673 1698}%
\special{fp}%
%
\special{pn 8}%
\special{pa 4569 1567}%
\special{pa 4575 1698}%
\special{fp}%
%
\special{pn 8}%
\special{pa 4490 1567}%
\special{pa 4490 1687}%
\special{fp}%
%
\special{pn 8}%
\special{pa 4402 1567}%
\special{pa 4402 2332}%
\special{fp}%
%
\special{pn 8}%
\special{pa 4351 1567}%
\special{pa 4351 2346}%
\special{fp}%
%
\special{pn 8}%
\special{pa 4496 2234}%
\special{pa 4502 2346}%
\special{fp}%
%
\special{pn 8}%
\special{pa 4607 2234}%
\special{pa 4613 2346}%
\special{fp}%
%
\special{pn 8}%
\special{pa 4751 2234}%
\special{pa 4751 2346}%
\special{fp}%
%
\special{pn 8}%
\special{pa 4686 2240}%
\special{pa 4686 2332}%
\special{fp}%
%
\special{pn 8}%
\special{pa 4904 2496}%
\special{pa 4772 2496}%
\special{fp}%
%
\special{pn 8}%
\special{pa 4857 2503}%
\special{pa 4865 2674}%
\special{fp}%
%
\special{pn 8}%
\special{pa 1543 642}%
\special{pa 3384 642}%
\special{pa 3384 2700}%
\special{pa 1543 2700}%
\special{pa 1543 642}%
\special{dt 0.045}%
%
\special{pn 8}%
\special{pa 1794 1040}%
\special{pa 3125 1060}%
\special{pa 3125 2510}%
\special{pa 2437 2528}%
\special{pa 2432 2112}%
\special{pa 2631 2112}%
\special{pa 2637 1965}%
\special{pa 1794 1965}%
\special{pa 1787 1956}%
\special{pa 1794 1060}%
\special{pa 1794 1060}%
\special{pa 1794 1060}%
\special{fp}%
%
\special{pn 8}%
\special{pa -696 624}%
\special{pa 1144 624}%
\special{pa 1144 2682}%
\special{pa -696 2682}%
\special{pa -696 624}%
\special{dt 0.045}%
%
\special{pn 8}%
\special{pa -447 1024}%
\special{pa 886 1043}%
\special{pa 886 2492}%
\special{pa 198 2510}%
\special{pa 191 2095}%
\special{pa 392 2095}%
\special{pa 398 1947}%
\special{pa -447 1947}%
\special{pa -452 1939}%
\special{pa -447 1043}%
\special{pa -447 1043}%
\special{pa -447 1043}%
\special{fp}%
%
\special{pn 8}%
\special{pa -273 1022}%
\special{pa -273 1956}%
\special{fp}%
%
\special{pn 8}%
\special{pa -79 1039}%
\special{pa -79 1905}%
\special{fp}%
%
\special{pn 8}%
\special{pa 103 1057}%
\special{pa 116 1905}%
\special{fp}%
%
\special{pn 8}%
\special{pa 309 1091}%
\special{pa 309 1905}%
\special{fp}%
%
\special{pn 8}%
\special{pa 539 1091}%
\special{pa 539 2458}%
\special{fp}%
%
\special{pn 8}%
\special{pa 721 1057}%
\special{pa 721 2441}%
\special{fp}%
%
\special{pn 8}%
\special{pa 333 2129}%
\special{pa 333 2441}%
\special{fp}%
\put(25.2000,-36.4000){\makebox(0,0)[lb]{Figure4.7.(1)}}%
\put(3.7000,-30.9000){\makebox(0,0)[lb]{t=-0.6}}%
\put(20.7000,-31.4000){\makebox(0,0)[lb]{t=-0.4}}%
\put(44.7000,-31.3000){\makebox(0,0)[lb]{t=-0.2}}%
\put(65.6000,-31.7000){\makebox(0,0)[lb]{t=0}}%
\end{picture}%
}

\np
\rotatebox[origin=c]{90}{
\unitlength 0.1in
\begin{picture}(62.36,28.32)(14.74,-34.60)
%
\special{pn 8}%
\special{pa 6899 2145}%
\special{pa 6899 2457}%
\special{fp}%
%
\special{pn 8}%
\special{pa 7287 1073}%
\special{pa 7287 2457}%
\special{fp}%
%
\special{pn 8}%
\special{pa 7105 1107}%
\special{pa 7105 2474}%
\special{fp}%
%
\special{pn 8}%
\special{pa 6875 1107}%
\special{pa 6875 1921}%
\special{fp}%
%
\special{pn 8}%
\special{pa 6669 1073}%
\special{pa 6682 1921}%
\special{fp}%
%
\special{pn 8}%
\special{pa 6487 1055}%
\special{pa 6487 1921}%
\special{fp}%
%
\special{pn 8}%
\special{pa 6293 1038}%
\special{pa 6293 1972}%
\special{fp}%
%
\special{pn 8}%
\special{pa 6119 1040}%
\special{pa 7452 1059}%
\special{pa 7452 2508}%
\special{pa 6764 2526}%
\special{pa 6757 2111}%
\special{pa 6958 2111}%
\special{pa 6964 1963}%
\special{pa 6119 1963}%
\special{pa 6114 1955}%
\special{pa 6119 1059}%
\special{pa 6119 1059}%
\special{pa 6119 1059}%
\special{fp}%
%
\special{pn 8}%
\special{pa 5870 640}%
\special{pa 7710 640}%
\special{pa 7710 2698}%
\special{pa 5870 2698}%
\special{pa 5870 640}%
\special{dt 0.045}%
%
\special{pn 8}%
\special{pa 3980 1030}%
\special{pa 5311 1050}%
\special{pa 5311 2500}%
\special{pa 4623 2518}%
\special{pa 4618 2102}%
\special{pa 4817 2102}%
\special{pa 4823 1955}%
\special{pa 3980 1955}%
\special{pa 3973 1946}%
\special{pa 3980 1050}%
\special{pa 3980 1050}%
\special{pa 3980 1050}%
\special{fp}%
%
\special{pn 8}%
\special{pa 2493 2485}%
\special{pa 2493 2656}%
\special{fp}%
%
\special{pn 8}%
\special{pa 2651 1945}%
\special{pa 1730 1945}%
\special{pa 1744 1028}%
\special{pa 3183 1028}%
\special{pa 3183 1945}%
\special{pa 3178 2471}%
\special{pa 2599 2471}%
\special{pa 2599 2656}%
\special{pa 2599 2656}%
\special{fp}%
%
\special{pn 8}%
\special{pa 1474 628}%
\special{pa 3467 628}%
\special{pa 3467 2663}%
\special{pa 1474 2663}%
\special{pa 1474 628}%
\special{dt 0.045}%
%
\special{pn 4}%
\special{sh 0}%
\special{pa 2362 995}%
\special{pa 2558 995}%
\special{pa 2558 1075}%
\special{pa 2362 1075}%
\special{pa 2362 995}%
\special{ip}%
%
\special{pn 4}%
\special{sh 0}%
\special{pa 1955 1920}%
\special{pa 2152 1920}%
\special{pa 2152 1997}%
\special{pa 1955 1997}%
\special{pa 1955 1920}%
\special{ip}%
%
\special{pn 8}%
\special{pa 2526 628}%
\special{pa 2526 1680}%
\special{pa 2126 1680}%
\special{pa 2126 2208}%
\special{pa 2440 2208}%
\special{pa 2447 2077}%
\special{pa 2658 2077}%
\special{pa 2651 1945}%
\special{pa 2651 1945}%
\special{pa 2651 1945}%
\special{fp}%
%
\special{pn 8}%
\special{pa 2388 628}%
\special{pa 2396 1543}%
\special{pa 1994 1543}%
\special{pa 1994 2335}%
\special{pa 2447 2335}%
\special{pa 2447 2656}%
\special{pa 2447 2656}%
\special{pa 2447 2656}%
\special{fp}%
%
\special{pn 8}%
\special{pa 2474 628}%
\special{pa 2474 1676}%
\special{fp}%
%
\special{pn 8}%
\special{pa 2428 628}%
\special{pa 2428 1669}%
\special{fp}%
%
\special{pn 8}%
\special{pa 2362 1543}%
\special{pa 2362 1680}%
\special{fp}%
%
\special{pn 8}%
\special{pa 2258 1549}%
\special{pa 2264 1680}%
\special{fp}%
%
\special{pn 8}%
\special{pa 2179 1549}%
\special{pa 2179 1669}%
\special{fp}%
%
\special{pn 8}%
\special{pa 2091 1549}%
\special{pa 2091 2314}%
\special{fp}%
%
\special{pn 8}%
\special{pa 2040 1549}%
\special{pa 2040 2328}%
\special{fp}%
%
\special{pn 8}%
\special{pa 2185 2216}%
\special{pa 2191 2328}%
\special{fp}%
%
\special{pn 8}%
\special{pa 2296 2216}%
\special{pa 2302 2328}%
\special{fp}%
%
\special{pn 8}%
\special{pa 2440 2216}%
\special{pa 2440 2328}%
\special{fp}%
%
\special{pn 8}%
\special{pa 2375 2222}%
\special{pa 2375 2314}%
\special{fp}%
%
\special{pn 8}%
\special{pa 2593 2478}%
\special{pa 2461 2478}%
\special{fp}%
%
\special{pn 8}%
\special{pa 2546 2485}%
\special{pa 2554 2656}%
\special{fp}%
%
\special{pn 8}%
\special{pa 3729 632}%
\special{pa 5570 632}%
\special{pa 5570 2690}%
\special{pa 3729 2690}%
\special{pa 3729 632}%
\special{dt 0.045}%
\put(34.8000,-36.3000){\makebox(0,0)[lb]{Figure 4.7.(2)}}%
\put(21.5000,-31.6000){\makebox(0,0)[lb]{t=0.2}}%
\put(42.9000,-31.8000){\makebox(0,0)[lb]{t=0.4}}%
\put(65.4000,-32.0000){\makebox(0,0)[lb]{t=0.6}}%
\end{picture}%
}

\np      
\rotatebox[origin=c]{90}{\input {4.8.1.tex}}
       
\np        
\rotatebox[origin=c]{90}{
\unitlength 0.1in
\begin{picture}(63.60,29.90)(7.40,-40.50)
%
\special{pn 8}%
\special{pa 3050 3550}%
\special{pa 5060 3550}%
\special{pa 5060 1060}%
\special{pa 3050 1060}%
\special{pa 3050 3550}%
\special{dt 0.045}%
%
\special{pn 8}%
\special{pa 1919 2824}%
\special{pa 998 2824}%
\special{pa 1012 1907}%
\special{pa 2451 1907}%
\special{pa 2451 2824}%
\special{pa 2446 3350}%
\special{pa 1867 3350}%
\special{pa 1867 3535}%
\special{pa 1867 3535}%
\special{fp}%
%
\special{pn 8}%
\special{pa 1861 3357}%
\special{pa 1729 3357}%
\special{fp}%
%
\special{pn 8}%
\special{pa 1643 3101}%
\special{pa 1643 3193}%
\special{fp}%
%
\special{pn 8}%
\special{pa 1708 3095}%
\special{pa 1708 3207}%
\special{fp}%
%
\special{pn 8}%
\special{pa 1564 3095}%
\special{pa 1570 3207}%
\special{fp}%
%
\special{pn 8}%
\special{pa 1453 3095}%
\special{pa 1459 3207}%
\special{fp}%
%
\special{pn 8}%
\special{pa 1447 2428}%
\special{pa 1447 2548}%
\special{fp}%
%
\special{pn 8}%
\special{pa 1526 2428}%
\special{pa 1532 2559}%
\special{fp}%
%
\special{pn 8}%
\special{pa 1630 2422}%
\special{pa 1630 2559}%
\special{fp}%
%
\special{pn 4}%
\special{sh 0}%
\special{pa 1630 1874}%
\special{pa 1826 1874}%
\special{pa 1826 1954}%
\special{pa 1630 1954}%
\special{pa 1630 1874}%
\special{ip}%
%
\special{pn 8}%
\special{pa 1761 3364}%
\special{pa 1761 3535}%
\special{fp}%
%
\special{pn 8}%
\special{pa 1814 3364}%
\special{pa 1822 3535}%
\special{fp}%
%
\special{pn 8}%
\special{pa 3411 1891}%
\special{pa 4742 1911}%
\special{pa 4742 3361}%
\special{pa 4054 3379}%
\special{pa 4049 2963}%
\special{pa 4248 2963}%
\special{pa 4254 2816}%
\special{pa 3411 2816}%
\special{pa 3404 2807}%
\special{pa 3411 1911}%
\special{pa 3411 1911}%
\special{pa 3411 1911}%
\special{fp}%
%
\special{pn 8}%
\special{pa 1794 1507}%
\special{pa 1794 2559}%
\special{pa 1394 2559}%
\special{pa 1394 3087}%
\special{pa 1708 3087}%
\special{pa 1715 2956}%
\special{pa 1926 2956}%
\special{pa 1919 2824}%
\special{pa 1919 2824}%
\special{pa 1919 2824}%
\special{fp}%
%
\special{pn 8}%
\special{pa 1656 1507}%
\special{pa 1664 2422}%
\special{pa 1262 2422}%
\special{pa 1262 3214}%
\special{pa 1715 3214}%
\special{pa 1715 3535}%
\special{pa 1715 3535}%
\special{pa 1715 3535}%
\special{fp}%
%
\special{pn 8}%
\special{pa 1359 2428}%
\special{pa 1359 3193}%
\special{fp}%
%
\special{pn 8}%
\special{pa 1308 2428}%
\special{pa 1308 3207}%
\special{fp}%
%
\special{pn 8}%
\special{pa 1742 1507}%
\special{pa 1742 2555}%
\special{fp}%
%
\special{pn 8}%
\special{pa 1696 1507}%
\special{pa 1696 2548}%
\special{fp}%
%
\special{pn 4}%
\special{sh 0}%
\special{pa 1223 2799}%
\special{pa 1420 2799}%
\special{pa 1420 2876}%
\special{pa 1223 2876}%
\special{pa 1223 2799}%
\special{ip}%
%
\special{pn 4}%
\special{pa 1170 2820}%
\special{pa 1660 2820}%
\special{ip}%
%
\special{pn 8}%
\special{pa 1130 2820}%
\special{pa 1700 2820}%
\special{fp}%
%
\special{pn 8}%
\special{pa 1790 1500}%
\special{pa 1650 1340}%
\special{fp}%
%
\special{pn 8}%
\special{pa 1650 1510}%
\special{pa 1690 1450}%
\special{fp}%
%
\special{pn 8}%
\special{pa 1750 1390}%
\special{pa 1800 1340}%
\special{fp}%
%
\special{pn 8}%
\special{pa 1800 1340}%
\special{pa 1650 1210}%
\special{fp}%
%
\special{pn 8}%
\special{pa 1660 1360}%
\special{pa 1710 1310}%
\special{fp}%
%
\special{pn 8}%
\special{pa 1750 1260}%
\special{pa 1800 1200}%
\special{fp}%
%
\special{pn 8}%
\special{pa 1790 1210}%
\special{pa 1790 1060}%
\special{fp}%
%
\special{pn 8}%
\special{pa 1660 1220}%
\special{pa 1660 1070}%
\special{fp}%
%
\special{pn 8}%
\special{pa 1710 1520}%
\special{pa 1700 1490}%
\special{fp}%
%
\special{pn 8}%
\special{pa 1690 1080}%
\special{pa 1690 1210}%
\special{fp}%
%
\special{pn 8}%
\special{pa 1750 1070}%
\special{pa 1750 1230}%
\special{fp}%
%
\special{pn 8}%
\special{pa 1730 1330}%
\special{pa 1730 1370}%
\special{fp}%
%
\special{pn 8}%
\special{pa 740 3550}%
\special{pa 2750 3550}%
\special{pa 2750 1060}%
\special{pa 740 1060}%
\special{pa 740 3550}%
\special{dt 0.045}%
%
\special{pn 8}%
\special{pa 6768 1918}%
\special{pa 6768 3302}%
\special{fp}%
%
\special{pn 8}%
\special{pa 6586 1952}%
\special{pa 6586 3319}%
\special{fp}%
%
\special{pn 8}%
\special{pa 6356 1952}%
\special{pa 6356 2766}%
\special{fp}%
%
\special{pn 8}%
\special{pa 5600 1885}%
\special{pa 6933 1904}%
\special{pa 6933 3353}%
\special{pa 6245 3371}%
\special{pa 6238 2956}%
\special{pa 6439 2956}%
\special{pa 6445 2808}%
\special{pa 5600 2808}%
\special{pa 5595 2800}%
\special{pa 5600 1904}%
\special{pa 5600 1904}%
\special{pa 5600 1904}%
\special{fp}%
%
\special{pn 8}%
\special{pa 6380 2990}%
\special{pa 6380 3302}%
\special{fp}%
%
\special{pn 8}%
\special{pa 5390 3570}%
\special{pa 7100 3570}%
\special{pa 7100 1060}%
\special{pa 5390 1060}%
\special{pa 5390 3570}%
\special{dt 0.045}%
%
\special{pn 8}%
\special{pa 6150 1918}%
\special{pa 6163 2766}%
\special{fp}%
%
\special{pn 8}%
\special{pa 5944 1907}%
\special{pa 5944 2773}%
\special{fp}%
%
\special{pn 8}%
\special{pa 5750 1890}%
\special{pa 5750 2824}%
\special{fp}%
\put(32.2000,-42.2000){\makebox(0,0)[lb]{Figure 4.8.(2)}}%
\put(13.5000,-39.1000){\makebox(0,0)[lb]{t=0.2}}%
\put(34.2000,-39.0000){\makebox(0,0)[lb]{t=0.4}}%
\put(59.2000,-39.1000){\makebox(0,0)[lb]{t=0.6}}%
\end{picture}%
}

\np
\rotatebox[origin=c]{90}{
\unitlength 0.1in
\begin{picture}(85.90,29.60)(-3.20,-40.90)
%
\special{pn 8}%
\special{pa 4290 1990}%
\special{pa 5740 1980}%
\special{pa 5730 3430}%
\special{pa 5150 3430}%
\special{pa 5150 3620}%
\special{pa 5150 3620}%
\special{fp}%
%
\special{pn 8}%
\special{pa 5090 1420}%
\special{pa 5120 1450}%
\special{fp}%
%
\special{pn 8}%
\special{pa 4950 1300}%
\special{pa 4950 1150}%
\special{fp}%
%
\special{pn 8}%
\special{pa 5080 1290}%
\special{pa 5080 1140}%
\special{fp}%
%
\special{pn 8}%
\special{pa 5040 1340}%
\special{pa 5090 1280}%
\special{fp}%
%
\special{pn 8}%
\special{pa 4950 1440}%
\special{pa 5000 1390}%
\special{fp}%
%
\special{pn 8}%
\special{pa 5090 1420}%
\special{pa 4940 1290}%
\special{fp}%
%
\special{pn 8}%
\special{pa 4990 3160}%
\special{pa 4990 3030}%
\special{pa 5210 3040}%
\special{pa 5190 2900}%
\special{pa 4280 2890}%
\special{pa 4280 1980}%
\special{pa 4280 1980}%
\special{fp}%
%
\special{pn 8}%
\special{pa 5200 2900}%
\special{pa 5210 3040}%
\special{pa 5000 3030}%
\special{pa 5000 3160}%
\special{pa 5000 3160}%
\special{fp}%
%
\special{pn 4}%
\special{pa 5000 3300}%
\special{pa 5000 3630}%
\special{ip}%
%
\special{pn 8}%
\special{pa 5000 3300}%
\special{pa 5000 3630}%
\special{fp}%
%
\special{pn 8}%
\special{pa 4020 3630}%
\special{pa 6030 3630}%
\special{pa 6030 1140}%
\special{pa 4020 1140}%
\special{pa 4020 3630}%
\special{dt 0.045}%
%
\special{pn 8}%
\special{pa 1770 3630}%
\special{pa 3780 3630}%
\special{pa 3780 1140}%
\special{pa 1770 1140}%
\special{pa 1770 3630}%
\special{dt 0.045}%
%
\special{pn 8}%
\special{pa 5050 1150}%
\special{pa 5050 1310}%
\special{fp}%
%
\special{pn 8}%
\special{pa 4990 1150}%
\special{pa 4990 1280}%
\special{fp}%
%
\special{pn 8}%
\special{pa 1058 1988}%
\special{pa 1058 3372}%
\special{fp}%
%
\special{pn 8}%
\special{pa 876 2022}%
\special{pa 876 3389}%
\special{fp}%
%
\special{pn 8}%
\special{pa 646 2022}%
\special{pa 646 2836}%
\special{fp}%
%
\special{pn 8}%
\special{pa -110 1955}%
\special{pa 1223 1974}%
\special{pa 1223 3423}%
\special{pa 535 3441}%
\special{pa 528 3026}%
\special{pa 729 3026}%
\special{pa 735 2878}%
\special{pa -110 2878}%
\special{pa -115 2870}%
\special{pa -110 1974}%
\special{pa -110 1974}%
\special{pa -110 1974}%
\special{fp}%
%
\special{pn 8}%
\special{pa 2131 1971}%
\special{pa 3462 1991}%
\special{pa 3462 3441}%
\special{pa 2774 3459}%
\special{pa 2769 3043}%
\special{pa 2968 3043}%
\special{pa 2974 2896}%
\special{pa 2131 2896}%
\special{pa 2124 2887}%
\special{pa 2131 1991}%
\special{pa 2131 1991}%
\special{pa 2131 1991}%
\special{fp}%
%
\special{pn 8}%
\special{pa 670 3060}%
\special{pa 670 3372}%
\special{fp}%
%
\special{pn 8}%
\special{pa -320 3640}%
\special{pa 1390 3640}%
\special{pa 1390 1130}%
\special{pa -320 1130}%
\special{pa -320 3640}%
\special{dt 0.045}%
%
\special{pn 8}%
\special{pa 440 1988}%
\special{pa 453 2836}%
\special{fp}%
%
\special{pn 8}%
\special{pa 234 1977}%
\special{pa 234 2843}%
\special{fp}%
%
\special{pn 8}%
\special{pa 40 1960}%
\special{pa 40 2894}%
\special{fp}%
%
\special{pn 8}%
\special{pa 6530 1990}%
\special{pa 7980 1980}%
\special{pa 7970 3430}%
\special{pa 7390 3430}%
\special{pa 7390 3620}%
\special{pa 7390 3620}%
\special{fp}%
%
\special{pn 8}%
\special{pa 6260 3630}%
\special{pa 8270 3630}%
\special{pa 8270 1140}%
\special{pa 6260 1140}%
\special{pa 6260 3630}%
\special{dt 0.045}%
%
\special{pn 8}%
\special{pa 7240 3300}%
\special{pa 7240 3630}%
\special{fp}%
%
\special{pn 4}%
\special{pa 7240 3300}%
\special{pa 7240 3630}%
\special{ip}%
%
\special{pn 8}%
\special{pa 7440 2900}%
\special{pa 7450 3040}%
\special{pa 7240 3030}%
\special{pa 7240 3160}%
\special{pa 7240 3160}%
\special{fp}%
%
\special{pn 8}%
\special{pa 7230 3160}%
\special{pa 7230 3030}%
\special{pa 7450 3040}%
\special{pa 7430 2900}%
\special{pa 6520 2890}%
\special{pa 6520 1980}%
\special{pa 6520 1980}%
\special{fp}%
%
\special{pn 8}%
\special{pa 7330 1420}%
\special{pa 7180 1290}%
\special{fp}%
%
\special{pn 8}%
\special{pa 7190 1440}%
\special{pa 7240 1390}%
\special{fp}%
%
\special{pn 8}%
\special{pa 7320 1330}%
\special{pa 7370 1270}%
\special{fp}%
%
\special{pn 8}%
\special{pa 7360 1290}%
\special{pa 7360 1140}%
\special{fp}%
%
\special{pn 8}%
\special{pa 7190 1300}%
\special{pa 7190 1150}%
\special{fp}%
%
\special{pn 4}%
\special{sh 0}%
\special{pa 7130 1950}%
\special{pa 7230 1950}%
\special{pa 7230 2070}%
\special{pa 7130 2070}%
\special{pa 7130 1950}%
\special{ip}%
%
\special{pn 4}%
\special{sh 0}%
\special{pa 7300 1930}%
\special{pa 7400 1930}%
\special{pa 7400 2050}%
\special{pa 7300 2050}%
\special{pa 7300 1930}%
\special{ip}%
%
\special{pn 8}%
\special{pa 7190 1450}%
\special{pa 7180 2640}%
\special{pa 7660 2640}%
\special{pa 7660 3160}%
\special{pa 7230 3160}%
\special{pa 7230 3160}%
\special{fp}%
%
\special{pn 8}%
\special{pa 7240 3300}%
\special{pa 7820 3300}%
\special{pa 7820 2500}%
\special{pa 7360 2500}%
\special{pa 7360 1440}%
\special{pa 7360 1440}%
\special{fp}%
%
\special{pn 8}%
\special{pa 7330 1420}%
\special{pa 7360 1450}%
\special{fp}%
%
\special{pn 4}%
\special{sh 0}%
\special{pa 4910 2020}%
\special{pa 5160 2020}%
\special{pa 5160 1940}%
\special{pa 4910 1940}%
\special{pa 4910 2020}%
\special{ip}%
%
\special{pn 8}%
\special{pa 5000 3300}%
\special{pa 5580 3300}%
\special{pa 5580 2500}%
\special{pa 5120 2500}%
\special{pa 5120 1440}%
\special{pa 5120 1440}%
\special{fp}%
%
\special{pn 8}%
\special{pa 4950 1450}%
\special{pa 4940 2640}%
\special{pa 5420 2640}%
\special{pa 5420 3160}%
\special{pa 4990 3160}%
\special{pa 4990 3160}%
\special{fp}%
%
\special{pn 8}%
\special{pa 5000 1440}%
\special{pa 5000 2620}%
\special{fp}%
%
\special{pn 8}%
\special{pa 5080 1440}%
\special{pa 5070 2610}%
\special{fp}%
%
\special{pn 8}%
\special{pa 5180 2520}%
\special{pa 5180 2630}%
\special{fp}%
%
\special{pn 8}%
\special{pa 5310 2510}%
\special{pa 5310 2630}%
\special{fp}%
%
\special{pn 8}%
\special{pa 5490 2510}%
\special{pa 5490 3280}%
\special{fp}%
%
\special{pn 8}%
\special{pa 5000 3170}%
\special{pa 5000 3300}%
\special{fp}%
%
\special{pn 8}%
\special{pa 5350 3160}%
\special{pa 5350 3310}%
\special{fp}%
%
\special{pn 8}%
\special{pa 5170 3170}%
\special{pa 5170 3300}%
\special{fp}%
%
\special{pn 8}%
\special{pa 5080 3180}%
\special{pa 5080 3290}%
\special{fp}%
%
\special{pn 8}%
\special{pa 5270 3170}%
\special{pa 5270 3290}%
\special{fp}%
\put(31.2000,-42.6000){\makebox(0,0)[lb]{Figure 4.9.(1)}}%
\put(2.3000,-39.5000){\makebox(0,0)[lb]{t=-0.6}}%
\put(24.0000,-40.0000){\makebox(0,0)[lb]{t=-0.4}}%
\put(46.0000,-40.0000){\makebox(0,0)[lb]{t=-0.2}}%
\put(66.0000,-40.0000){\makebox(0,0)[lb]{t=0}}%
\end{picture}%
}

\np
\rotatebox[origin=c]{90}{
\unitlength 0.1in
\begin{picture}(65.80,31.70)(5.70,-41.10)
%
\special{pn 8}%
\special{pa 840 1790}%
\special{pa 2290 1780}%
\special{pa 2280 3230}%
\special{pa 1700 3230}%
\special{pa 1700 3420}%
\special{pa 1700 3420}%
\special{fp}%
%
\special{pn 8}%
\special{pa 1640 1220}%
\special{pa 1670 1250}%
\special{fp}%
%
\special{pn 8}%
\special{pa 1500 1100}%
\special{pa 1500 950}%
\special{fp}%
%
\special{pn 8}%
\special{pa 1630 1090}%
\special{pa 1630 940}%
\special{fp}%
%
\special{pn 8}%
\special{pa 1590 1140}%
\special{pa 1640 1080}%
\special{fp}%
%
\special{pn 8}%
\special{pa 1500 1240}%
\special{pa 1550 1190}%
\special{fp}%
%
\special{pn 8}%
\special{pa 1640 1220}%
\special{pa 1490 1090}%
\special{fp}%
%
\special{pn 8}%
\special{pa 1540 2960}%
\special{pa 1540 2830}%
\special{pa 1760 2840}%
\special{pa 1740 2700}%
\special{pa 830 2690}%
\special{pa 830 1780}%
\special{pa 830 1780}%
\special{fp}%
%
\special{pn 8}%
\special{pa 1750 2700}%
\special{pa 1760 2840}%
\special{pa 1550 2830}%
\special{pa 1550 2960}%
\special{pa 1550 2960}%
\special{fp}%
%
\special{pn 4}%
\special{pa 1550 3100}%
\special{pa 1550 3430}%
\special{ip}%
%
\special{pn 8}%
\special{pa 1550 3100}%
\special{pa 1550 3430}%
\special{fp}%
%
\special{pn 8}%
\special{pa 570 3430}%
\special{pa 2580 3430}%
\special{pa 2580 940}%
\special{pa 570 940}%
\special{pa 570 3430}%
\special{dt 0.045}%
%
\special{pn 8}%
\special{pa 3020 3430}%
\special{pa 5030 3430}%
\special{pa 5030 940}%
\special{pa 3020 940}%
\special{pa 3020 3430}%
\special{dt 0.045}%
%
\special{pn 8}%
\special{pa 1600 950}%
\special{pa 1600 1110}%
\special{fp}%
%
\special{pn 8}%
\special{pa 1540 950}%
\special{pa 1540 1080}%
\special{fp}%
%
\special{pn 8}%
\special{pa 6818 1798}%
\special{pa 6818 3182}%
\special{fp}%
%
\special{pn 8}%
\special{pa 6636 1832}%
\special{pa 6636 3199}%
\special{fp}%
%
\special{pn 8}%
\special{pa 6406 1832}%
\special{pa 6406 2646}%
\special{fp}%
%
\special{pn 8}%
\special{pa 5650 1765}%
\special{pa 6984 1784}%
\special{pa 6984 3233}%
\special{pa 6296 3251}%
\special{pa 6288 2836}%
\special{pa 6488 2836}%
\special{pa 6496 2688}%
\special{pa 5650 2688}%
\special{pa 5644 2680}%
\special{pa 5650 1784}%
\special{pa 5650 1784}%
\special{pa 5650 1784}%
\special{fp}%
%
\special{pn 8}%
\special{pa 3382 1771}%
\special{pa 4712 1791}%
\special{pa 4712 3241}%
\special{pa 4024 3259}%
\special{pa 4018 2843}%
\special{pa 4218 2843}%
\special{pa 4224 2696}%
\special{pa 3382 2696}%
\special{pa 3374 2687}%
\special{pa 3382 1791}%
\special{pa 3382 1791}%
\special{pa 3382 1791}%
\special{fp}%
%
\special{pn 8}%
\special{pa 6430 2870}%
\special{pa 6430 3182}%
\special{fp}%
%
\special{pn 8}%
\special{pa 5440 3450}%
\special{pa 7150 3450}%
\special{pa 7150 940}%
\special{pa 5440 940}%
\special{pa 5440 3450}%
\special{dt 0.045}%
%
\special{pn 8}%
\special{pa 6200 1798}%
\special{pa 6212 2646}%
\special{fp}%
%
\special{pn 8}%
\special{pa 5994 1787}%
\special{pa 5994 2653}%
\special{fp}%
%
\special{pn 8}%
\special{pa 5800 1770}%
\special{pa 5800 2704}%
\special{fp}%
%
\special{pn 4}%
\special{sh 0}%
\special{pa 1460 1820}%
\special{pa 1710 1820}%
\special{pa 1710 1740}%
\special{pa 1460 1740}%
\special{pa 1460 1820}%
\special{ip}%
%
\special{pn 8}%
\special{pa 1550 3100}%
\special{pa 2130 3100}%
\special{pa 2130 2300}%
\special{pa 1670 2300}%
\special{pa 1670 1240}%
\special{pa 1670 1240}%
\special{fp}%
%
\special{pn 8}%
\special{pa 1500 1250}%
\special{pa 1490 2440}%
\special{pa 1970 2440}%
\special{pa 1970 2960}%
\special{pa 1540 2960}%
\special{pa 1540 2960}%
\special{fp}%
%
\special{pn 8}%
\special{pa 1550 1240}%
\special{pa 1550 2420}%
\special{fp}%
%
\special{pn 8}%
\special{pa 1630 1240}%
\special{pa 1620 2410}%
\special{fp}%
%
\special{pn 8}%
\special{pa 1730 2320}%
\special{pa 1730 2430}%
\special{fp}%
%
\special{pn 8}%
\special{pa 1860 2310}%
\special{pa 1860 2430}%
\special{fp}%
%
\special{pn 8}%
\special{pa 2040 2310}%
\special{pa 2040 3080}%
\special{fp}%
%
\special{pn 8}%
\special{pa 1550 2970}%
\special{pa 1550 3100}%
\special{fp}%
%
\special{pn 8}%
\special{pa 1900 2960}%
\special{pa 1900 3110}%
\special{fp}%
%
\special{pn 8}%
\special{pa 1720 2970}%
\special{pa 1720 3100}%
\special{fp}%
%
\special{pn 8}%
\special{pa 1630 2980}%
\special{pa 1630 3090}%
\special{fp}%
%
\special{pn 8}%
\special{pa 1820 2970}%
\special{pa 1820 3090}%
\special{fp}%
\put(29.2000,-42.8000){\makebox(0,0)[lb]{Figure 4.9.(2)}}%
\put(12.0000,-38.0000){\makebox(0,0)[lb]{t=0.2}}%
\put(36.0000,-38.0000){\makebox(0,0)[lb]{t=0.4}}%
\put(60.0000,-38.0000){\makebox(0,0)[lb]{t=0.6}}%
\end{picture}%
}

\np
\unitlength 0.1in
\begin{picture}(47.30,25.50)(9.40,-27.50)
%
\special{pn 8}%
\special{pa 3200 210}%
\special{pa 3200 1430}%
\special{fp}%
%
\special{pn 8}%
\special{pa 3470 200}%
\special{pa 3470 1420}%
\special{fp}%
%
\special{pn 8}%
\special{pa 2820 200}%
\special{pa 3900 200}%
\special{pa 3900 1420}%
\special{pa 2820 1420}%
\special{pa 2820 200}%
\special{dt 0.045}%
%
\special{pn 8}%
\special{pa 1320 210}%
\special{pa 1320 1430}%
\special{fp}%
%
\special{pn 8}%
\special{pa 1590 200}%
\special{pa 1590 1420}%
\special{fp}%
%
\special{pn 8}%
\special{pa 940 200}%
\special{pa 2020 200}%
\special{pa 2020 1420}%
\special{pa 940 1420}%
\special{pa 940 200}%
\special{dt 0.045}%
%
\special{pn 8}%
\special{pa 1370 210}%
\special{pa 1370 1420}%
\special{fp}%
%
\special{pn 8}%
\special{pa 1470 210}%
\special{pa 1460 1420}%
\special{fp}%
%
\special{pn 8}%
\special{pa 1530 210}%
\special{pa 1540 1400}%
\special{fp}%
%
\special{pn 8}%
\special{pa 4970 210}%
\special{pa 4970 1430}%
\special{fp}%
%
\special{pn 8}%
\special{pa 5240 200}%
\special{pa 5240 1420}%
\special{fp}%
%
\special{pn 8}%
\special{pa 4590 200}%
\special{pa 5670 200}%
\special{pa 5670 1420}%
\special{pa 4590 1420}%
\special{pa 4590 200}%
\special{dt 0.045}%
%
\special{pn 8}%
\special{pa 5020 210}%
\special{pa 5020 1420}%
\special{fp}%
%
\special{pn 8}%
\special{pa 5120 210}%
\special{pa 5110 1420}%
\special{fp}%
%
\special{pn 8}%
\special{pa 5180 210}%
\special{pa 5190 1400}%
\special{fp}%
\put(27.0000,-29.2000){\makebox(0,0)[lb]{Figure 4.10}}%
\put(12.0000,-20.0000){\makebox(0,0)[lb]{t=-0.5}}%
\put(32.0000,-20.0000){\makebox(0,0)[lb]{t=0}}%
\put(49.1000,-20.2000){\makebox(0,0)[lb]{t=0.5}}%
\end{picture}%

\np
\unitlength 0.1in
\begin{picture}(47.10,23.10)(10.00,-26.40)
%
\special{pn 8}%
\special{pa 3270 340}%
\special{pa 3270 540}%
\special{fp}%
%
\special{pn 8}%
\special{pa 3270 560}%
\special{pa 3530 830}%
\special{fp}%
%
\special{pn 8}%
\special{pa 3530 840}%
\special{pa 3460 930}%
\special{fp}%
%
\special{pn 8}%
\special{pa 3520 340}%
\special{pa 3520 520}%
\special{fp}%
%
\special{pn 8}%
\special{pa 3520 530}%
\special{pa 3440 630}%
\special{fp}%
%
\special{pn 8}%
\special{pa 3360 760}%
\special{pa 3260 880}%
\special{fp}%
%
\special{pn 8}%
\special{pa 3280 880}%
\special{pa 3520 1110}%
\special{fp}%
%
\special{pn 8}%
\special{pa 3380 1060}%
\special{pa 3260 1180}%
\special{fp}%
%
\special{pn 8}%
\special{pa 3280 1180}%
\special{pa 3280 1510}%
\special{fp}%
%
\special{pn 8}%
\special{pa 3530 1110}%
\special{pa 3530 1500}%
\special{fp}%
%
\special{pn 8}%
\special{pa 2920 330}%
\special{pa 3980 330}%
\special{pa 3980 1540}%
\special{pa 2920 1540}%
\special{pa 2920 330}%
\special{dt 0.045}%
%
\special{pn 8}%
\special{pa 1350 340}%
\special{pa 1350 540}%
\special{fp}%
%
\special{pn 8}%
\special{pa 1350 560}%
\special{pa 1610 830}%
\special{fp}%
%
\special{pn 8}%
\special{pa 1610 840}%
\special{pa 1540 930}%
\special{fp}%
%
\special{pn 8}%
\special{pa 1600 340}%
\special{pa 1600 520}%
\special{fp}%
%
\special{pn 8}%
\special{pa 1600 530}%
\special{pa 1520 630}%
\special{fp}%
%
\special{pn 8}%
\special{pa 1440 760}%
\special{pa 1340 880}%
\special{fp}%
%
\special{pn 8}%
\special{pa 1360 880}%
\special{pa 1600 1110}%
\special{fp}%
%
\special{pn 8}%
\special{pa 1460 1060}%
\special{pa 1340 1180}%
\special{fp}%
%
\special{pn 8}%
\special{pa 1360 1180}%
\special{pa 1360 1510}%
\special{fp}%
%
\special{pn 8}%
\special{pa 1610 1110}%
\special{pa 1610 1500}%
\special{fp}%
%
\special{pn 8}%
\special{pa 1000 330}%
\special{pa 2060 330}%
\special{pa 2060 1540}%
\special{pa 1000 1540}%
\special{pa 1000 330}%
\special{dt 0.045}%
%
\special{pn 8}%
\special{pa 1440 340}%
\special{pa 1440 600}%
\special{fp}%
%
\special{pn 8}%
\special{pa 1530 330}%
\special{pa 1530 560}%
\special{fp}%
%
\special{pn 8}%
\special{pa 1430 1130}%
\special{pa 1440 1540}%
\special{fp}%
%
\special{pn 8}%
\special{pa 1550 1110}%
\special{pa 1540 1540}%
\special{fp}%
%
\special{pn 8}%
\special{pa 1440 840}%
\special{pa 1440 890}%
\special{fp}%
%
\special{pn 8}%
\special{pa 1530 800}%
\special{pa 1530 870}%
\special{fp}%
%
\special{pn 8}%
\special{pa 5000 340}%
\special{pa 5000 540}%
\special{fp}%
%
\special{pn 8}%
\special{pa 5000 560}%
\special{pa 5260 830}%
\special{fp}%
%
\special{pn 8}%
\special{pa 5260 840}%
\special{pa 5190 930}%
\special{fp}%
%
\special{pn 8}%
\special{pa 5250 340}%
\special{pa 5250 520}%
\special{fp}%
%
\special{pn 8}%
\special{pa 5250 530}%
\special{pa 5170 630}%
\special{fp}%
%
\special{pn 8}%
\special{pa 5090 760}%
\special{pa 4990 880}%
\special{fp}%
%
\special{pn 8}%
\special{pa 5010 880}%
\special{pa 5250 1110}%
\special{fp}%
%
\special{pn 8}%
\special{pa 5110 1060}%
\special{pa 4990 1180}%
\special{fp}%
%
\special{pn 8}%
\special{pa 5010 1180}%
\special{pa 5010 1510}%
\special{fp}%
%
\special{pn 8}%
\special{pa 5260 1110}%
\special{pa 5260 1500}%
\special{fp}%
%
\special{pn 8}%
\special{pa 4650 330}%
\special{pa 5710 330}%
\special{pa 5710 1540}%
\special{pa 4650 1540}%
\special{pa 4650 330}%
\special{dt 0.045}%
%
\special{pn 8}%
\special{pa 5090 340}%
\special{pa 5090 600}%
\special{fp}%
%
\special{pn 8}%
\special{pa 5180 330}%
\special{pa 5180 560}%
\special{fp}%
%
\special{pn 8}%
\special{pa 5080 1130}%
\special{pa 5090 1540}%
\special{fp}%
%
\special{pn 8}%
\special{pa 5200 1110}%
\special{pa 5190 1540}%
\special{fp}%
%
\special{pn 8}%
\special{pa 5090 840}%
\special{pa 5090 890}%
\special{fp}%
%
\special{pn 8}%
\special{pa 5180 800}%
\special{pa 5180 870}%
\special{fp}%
\put(29.8000,-28.1000){\makebox(0,0)[lb]{Figure 4.11}}%
\put(12.0000,-20.0000){\makebox(0,0)[lb]{t=-0.5}}%
\put(32.0000,-20.0000){\makebox(0,0)[lb]{t=0}}%
\put(49.1000,-20.2000){\makebox(0,0)[lb]{t=0.5}}%
\end{picture}%

\np
\input 4.12.tex

\np
\input 4.13.tex

\np
\unitlength 0.1in
\begin{picture}(56.10,45.70)(8.50,-46.30)
%
\special{pn 8}%
\special{ar 3510 320 560 250  0.0000000 6.2831853}%
%
\special{pn 20}%
\special{pa 3660 3950}%
\special{pa 3660 2730}%
\special{fp}%
%
\special{pn 20}%
\special{pa 3340 3970}%
\special{pa 3340 2740}%
\special{fp}%
%
\special{pn 20}%
\special{pa 3350 4000}%
\special{pa 3360 3970}%
\special{pa 3384 3949}%
\special{pa 3413 3935}%
\special{pa 3444 3926}%
\special{pa 3475 3921}%
\special{pa 3507 3920}%
\special{pa 3539 3923}%
\special{pa 3570 3930}%
\special{pa 3601 3940}%
\special{pa 3627 3958}%
\special{pa 3647 3983}%
\special{pa 3648 4014}%
\special{pa 3630 4040}%
\special{pa 3603 4058}%
\special{pa 3573 4070}%
\special{pa 3542 4077}%
\special{pa 3510 4080}%
\special{pa 3478 4079}%
\special{pa 3447 4075}%
\special{pa 3416 4066}%
\special{pa 3387 4053}%
\special{pa 3362 4032}%
\special{pa 3350 4003}%
\special{pa 3350 4000}%
\special{sp}%
%
\special{pn 20}%
\special{pa 3340 360}%
\special{pa 3340 1580}%
\special{fp}%
%
\special{pn 20}%
\special{pa 3660 340}%
\special{pa 3660 1570}%
\special{fp}%
%
\special{pn 20}%
\special{ar 3500 310 150 80  0.0000000 6.2831853}%
%
\special{pn 8}%
\special{ar 3510 4080 560 250  0.0000000 6.2831853}%
%
\special{pn 8}%
\special{pa 4080 330}%
\special{pa 4080 4080}%
\special{fp}%
%
\special{pn 8}%
\special{pa 2950 340}%
\special{pa 2950 4080}%
\special{fp}%
%
\special{pn 8}%
\special{ar 1410 330 560 250  0.0000000 6.2831853}%
%
\special{pn 8}%
\special{ar 1410 4090 560 250  0.0000000 6.2831853}%
%
\special{pn 8}%
\special{pa 1980 340}%
\special{pa 1980 4090}%
\special{fp}%
%
\special{pn 8}%
\special{pa 850 350}%
\special{pa 850 4090}%
\special{fp}%
%
\special{pn 8}%
\special{ar 5890 310 560 250  0.0000000 6.2831853}%
%
\special{pn 8}%
\special{ar 5890 4070 560 250  0.0000000 6.2831853}%
%
\special{pn 8}%
\special{pa 6460 320}%
\special{pa 6460 4070}%
\special{fp}%
%
\special{pn 8}%
\special{pa 5330 330}%
\special{pa 5330 4070}%
\special{fp}%
%
\special{pn 20}%
\special{ar 3500 1590 150 80  0.0000000 6.2831853}%
%
\special{pn 20}%
\special{ar 3500 2710 150 80  0.0000000 6.2831853}%
%
\special{pn 20}%
\special{ar 3520 2070 560 250  0.0000000 6.2831853}%
%
\special{pn 8}%
\special{pa 3660 1600}%
\special{pa 5730 1600}%
\special{dt 0.045}%
\special{pa 5730 1600}%
\special{pa 5729 1600}%
\special{dt 0.045}%
%
\special{pn 8}%
\special{pa 3680 2720}%
\special{pa 5750 2720}%
\special{dt 0.045}%
\special{pa 5750 2720}%
\special{pa 5749 2720}%
\special{dt 0.045}%
%
\special{pn 20}%
\special{ar 5850 1600 150 80  0.0000000 6.2831853}%
%
\special{pn 20}%
\special{ar 5850 2730 150 80  0.0000000 6.2831853}%
%
\special{pn 20}%
\special{pa 5690 1610}%
\special{pa 5690 2700}%
\special{fp}%
%
\special{pn 20}%
\special{pa 6020 1640}%
\special{pa 6020 2730}%
\special{fp}%
\put(12.0000,-46.0000){\makebox(0,0)[lb]{t=-0.5}}%
\put(32.0000,-46.0000){\makebox(0,0)[lb]{t=0}}%
\put(56.0000,-46.0000){\makebox(0,0)[lb]{t=0.5}}%
\put(34.7000,-48.0000){\makebox(0,0)[lb]{Figure 9.1}}%
%
\special{pn 8}%
\special{ar 3510 320 560 250  0.0000000 6.2831853}%
%
\special{pn 8}%
\special{pa 2950 340}%
\special{pa 2950 4080}%
\special{fp}%
%
\special{pn 8}%
\special{pa 4080 330}%
\special{pa 4080 4080}%
\special{fp}%
%
\special{pn 8}%
\special{ar 3510 4080 560 250  0.0000000 6.2831853}%
%
\special{pn 20}%
\special{ar 3500 310 150 80  0.0000000 6.2831853}%
%
\special{pn 20}%
\special{pa 3660 340}%
\special{pa 3660 1570}%
\special{fp}%
%
\special{pn 20}%
\special{pa 3340 360}%
\special{pa 3340 1580}%
\special{fp}%
%
\special{pn 20}%
\special{pa 3350 4000}%
\special{pa 3360 3970}%
\special{pa 3384 3949}%
\special{pa 3413 3935}%
\special{pa 3444 3926}%
\special{pa 3475 3921}%
\special{pa 3507 3920}%
\special{pa 3539 3923}%
\special{pa 3570 3930}%
\special{pa 3601 3940}%
\special{pa 3627 3958}%
\special{pa 3647 3983}%
\special{pa 3648 4014}%
\special{pa 3630 4040}%
\special{pa 3603 4058}%
\special{pa 3573 4070}%
\special{pa 3542 4077}%
\special{pa 3510 4080}%
\special{pa 3478 4079}%
\special{pa 3447 4075}%
\special{pa 3416 4066}%
\special{pa 3387 4053}%
\special{pa 3362 4032}%
\special{pa 3350 4003}%
\special{pa 3350 4000}%
\special{sp}%
%
\special{pn 20}%
\special{pa 3340 3970}%
\special{pa 3340 2740}%
\special{fp}%
%
\special{pn 20}%
\special{pa 3660 3950}%
\special{pa 3660 2730}%
\special{fp}%
%
\special{pn 8}%
\special{ar 1410 330 560 250  0.0000000 6.2831853}%
%
\special{pn 8}%
\special{ar 1410 4090 560 250  0.0000000 6.2831853}%
%
\special{pn 8}%
\special{pa 1980 340}%
\special{pa 1980 4090}%
\special{fp}%
%
\special{pn 8}%
\special{pa 850 350}%
\special{pa 850 4090}%
\special{fp}%
%
\special{pn 8}%
\special{ar 5890 310 560 250  0.0000000 6.2831853}%
%
\special{pn 8}%
\special{ar 5890 4070 560 250  0.0000000 6.2831853}%
%
\special{pn 8}%
\special{pa 6460 320}%
\special{pa 6460 4070}%
\special{fp}%
%
\special{pn 8}%
\special{pa 5330 330}%
\special{pa 5330 4070}%
\special{fp}%
\end{picture}%

\np
\unitlength 0.1in
\begin{picture}(56.10,45.60)(8.50,-46.20)
%
\special{pn 8}%
\special{ar 3510 320 560 250  0.0000000 6.2831853}%
%
\special{pn 20}%
\special{pa 3660 3950}%
\special{pa 3660 2730}%
\special{fp}%
%
\special{pn 20}%
\special{pa 3340 3970}%
\special{pa 3340 2740}%
\special{fp}%
%
\special{pn 20}%
\special{pa 3350 4000}%
\special{pa 3360 3970}%
\special{pa 3384 3949}%
\special{pa 3413 3935}%
\special{pa 3444 3926}%
\special{pa 3475 3921}%
\special{pa 3507 3920}%
\special{pa 3539 3923}%
\special{pa 3570 3930}%
\special{pa 3601 3940}%
\special{pa 3627 3958}%
\special{pa 3647 3983}%
\special{pa 3648 4014}%
\special{pa 3630 4040}%
\special{pa 3603 4058}%
\special{pa 3573 4070}%
\special{pa 3542 4077}%
\special{pa 3510 4080}%
\special{pa 3478 4079}%
\special{pa 3447 4075}%
\special{pa 3416 4066}%
\special{pa 3387 4053}%
\special{pa 3362 4032}%
\special{pa 3350 4003}%
\special{pa 3350 4000}%
\special{sp}%
%
\special{pn 20}%
\special{pa 3340 360}%
\special{pa 3340 1580}%
\special{fp}%
%
\special{pn 20}%
\special{pa 3660 340}%
\special{pa 3660 1570}%
\special{fp}%
%
\special{pn 20}%
\special{ar 3500 310 150 80  0.0000000 6.2831853}%
%
\special{pn 8}%
\special{ar 3510 4080 560 250  0.0000000 6.2831853}%
%
\special{pn 8}%
\special{pa 4080 330}%
\special{pa 4080 4080}%
\special{fp}%
%
\special{pn 8}%
\special{pa 2950 340}%
\special{pa 2950 4080}%
\special{fp}%
%
\special{pn 8}%
\special{ar 1410 330 560 250  0.0000000 6.2831853}%
%
\special{pn 8}%
\special{ar 1410 4090 560 250  0.0000000 6.2831853}%
%
\special{pn 8}%
\special{pa 1980 340}%
\special{pa 1980 4090}%
\special{fp}%
%
\special{pn 8}%
\special{pa 850 350}%
\special{pa 850 4090}%
\special{fp}%
%
\special{pn 8}%
\special{ar 5890 310 560 250  0.0000000 6.2831853}%
%
\special{pn 8}%
\special{ar 5890 4070 560 250  0.0000000 6.2831853}%
%
\special{pn 8}%
\special{pa 6460 320}%
\special{pa 6460 4070}%
\special{fp}%
%
\special{pn 8}%
\special{pa 5330 330}%
\special{pa 5330 4070}%
\special{fp}%
%
\special{pn 20}%
\special{ar 3500 1590 150 80  0.0000000 6.2831853}%
%
\special{pn 20}%
\special{ar 3500 2710 150 80  0.0000000 6.2831853}%
%
\special{pn 20}%
\special{ar 3520 2070 560 250  0.0000000 6.2831853}%
\put(12.0000,-46.0000){\makebox(0,0)[lb]{t=-0.5}}%
\put(32.0000,-46.0000){\makebox(0,0)[lb]{t=0}}%
\put(56.0000,-46.0000){\makebox(0,0)[lb]{t=0.5}}%
%
\special{pn 8}%
\special{ar 3510 320 560 250  0.0000000 6.2831853}%
%
\special{pn 8}%
\special{pa 2950 340}%
\special{pa 2950 4080}%
\special{fp}%
%
\special{pn 8}%
\special{pa 4080 330}%
\special{pa 4080 4080}%
\special{fp}%
%
\special{pn 8}%
\special{ar 3510 4080 560 250  0.0000000 6.2831853}%
%
\special{pn 20}%
\special{ar 3500 310 150 80  0.0000000 6.2831853}%
%
\special{pn 20}%
\special{pa 3660 340}%
\special{pa 3660 1570}%
\special{fp}%
%
\special{pn 20}%
\special{pa 3340 360}%
\special{pa 3340 1580}%
\special{fp}%
%
\special{pn 20}%
\special{pa 3350 4000}%
\special{pa 3360 3970}%
\special{pa 3384 3949}%
\special{pa 3413 3935}%
\special{pa 3444 3926}%
\special{pa 3475 3921}%
\special{pa 3507 3920}%
\special{pa 3539 3923}%
\special{pa 3570 3930}%
\special{pa 3601 3940}%
\special{pa 3627 3958}%
\special{pa 3647 3983}%
\special{pa 3648 4014}%
\special{pa 3630 4040}%
\special{pa 3603 4058}%
\special{pa 3573 4070}%
\special{pa 3542 4077}%
\special{pa 3510 4080}%
\special{pa 3478 4079}%
\special{pa 3447 4075}%
\special{pa 3416 4066}%
\special{pa 3387 4053}%
\special{pa 3362 4032}%
\special{pa 3350 4003}%
\special{pa 3350 4000}%
\special{sp}%
%
\special{pn 20}%
\special{pa 3340 3970}%
\special{pa 3340 2740}%
\special{fp}%
%
\special{pn 20}%
\special{pa 3660 3950}%
\special{pa 3660 2730}%
\special{fp}%
%
\special{pn 8}%
\special{ar 1410 330 560 250  0.0000000 6.2831853}%
%
\special{pn 8}%
\special{ar 1410 4090 560 250  0.0000000 6.2831853}%
%
\special{pn 8}%
\special{pa 1980 340}%
\special{pa 1980 4090}%
\special{fp}%
%
\special{pn 8}%
\special{pa 850 350}%
\special{pa 850 4090}%
\special{fp}%
%
\special{pn 8}%
\special{ar 5890 310 560 250  0.0000000 6.2831853}%
%
\special{pn 8}%
\special{ar 5890 4070 560 250  0.0000000 6.2831853}%
%
\special{pn 8}%
\special{pa 6460 320}%
\special{pa 6460 4070}%
\special{fp}%
%
\special{pn 8}%
\special{pa 5330 330}%
\special{pa 5330 4070}%
\special{fp}%
%
\special{pn 8}%
\special{pa 3320 1590}%
\special{pa 1590 1590}%
\special{dt 0.045}%
\special{pa 1590 1590}%
\special{pa 1591 1590}%
\special{dt 0.045}%
%
\special{pn 8}%
\special{pa 3330 2740}%
\special{pa 1620 2740}%
\special{dt 0.045}%
\special{pa 1620 2740}%
\special{pa 1621 2740}%
\special{dt 0.045}%
%
\special{pn 20}%
\special{ar 1420 1590 150 80  0.0000000 6.2831853}%
%
\special{pn 20}%
\special{ar 1420 2720 150 80  0.0000000 6.2831853}%
%
\special{pn 20}%
\special{pa 1260 1600}%
\special{pa 1260 2690}%
\special{fp}%
%
\special{pn 20}%
\special{pa 1590 1630}%
\special{pa 1590 2720}%
\special{fp}%
\put(32.4000,-47.9000){\makebox(0,0)[lb]{Figure 9.2}}%
\end{picture}%

\np
\unitlength 0.1in
\begin{picture}(56.10,45.60)(8.50,-46.20)
%
\special{pn 8}%
\special{ar 3510 320 560 250  0.0000000 6.2831853}%
%
\special{pn 20}%
\special{pa 3350 4000}%
\special{pa 3360 3970}%
\special{pa 3384 3949}%
\special{pa 3413 3935}%
\special{pa 3444 3926}%
\special{pa 3475 3921}%
\special{pa 3507 3920}%
\special{pa 3539 3923}%
\special{pa 3570 3930}%
\special{pa 3601 3940}%
\special{pa 3627 3958}%
\special{pa 3647 3983}%
\special{pa 3648 4014}%
\special{pa 3630 4040}%
\special{pa 3603 4058}%
\special{pa 3573 4070}%
\special{pa 3542 4077}%
\special{pa 3510 4080}%
\special{pa 3478 4079}%
\special{pa 3447 4075}%
\special{pa 3416 4066}%
\special{pa 3387 4053}%
\special{pa 3362 4032}%
\special{pa 3350 4003}%
\special{pa 3350 4000}%
\special{sp}%
%
\special{pn 20}%
\special{ar 3500 310 150 80  0.0000000 6.2831853}%
%
\special{pn 8}%
\special{ar 3510 4080 560 250  0.0000000 6.2831853}%
%
\special{pn 8}%
\special{pa 4080 330}%
\special{pa 4080 4080}%
\special{fp}%
%
\special{pn 8}%
\special{pa 2950 340}%
\special{pa 2950 4080}%
\special{fp}%
%
\special{pn 8}%
\special{ar 1410 330 560 250  0.0000000 6.2831853}%
%
\special{pn 8}%
\special{ar 1410 4090 560 250  0.0000000 6.2831853}%
%
\special{pn 8}%
\special{pa 1980 340}%
\special{pa 1980 4090}%
\special{fp}%
%
\special{pn 8}%
\special{pa 850 350}%
\special{pa 850 4090}%
\special{fp}%
%
\special{pn 8}%
\special{ar 5890 310 560 250  0.0000000 6.2831853}%
%
\special{pn 8}%
\special{ar 5890 4070 560 250  0.0000000 6.2831853}%
%
\special{pn 8}%
\special{pa 6460 320}%
\special{pa 6460 4070}%
\special{fp}%
%
\special{pn 8}%
\special{pa 5330 330}%
\special{pa 5330 4070}%
\special{fp}%
%
\special{pn 20}%
\special{ar 3520 2070 560 250  0.0000000 6.2831853}%
\put(12.0000,-46.0000){\makebox(0,0)[lb]{t=-0.5}}%
\put(32.0000,-46.0000){\makebox(0,0)[lb]{t=0}}%
\put(56.0000,-46.0000){\makebox(0,0)[lb]{t=0.5}}%
%
\special{pn 8}%
\special{ar 3510 320 560 250  0.0000000 6.2831853}%
%
\special{pn 20}%
\special{pa 3350 4000}%
\special{pa 3360 3970}%
\special{pa 3384 3949}%
\special{pa 3413 3935}%
\special{pa 3444 3926}%
\special{pa 3475 3921}%
\special{pa 3507 3920}%
\special{pa 3539 3923}%
\special{pa 3570 3930}%
\special{pa 3601 3940}%
\special{pa 3627 3958}%
\special{pa 3647 3983}%
\special{pa 3648 4014}%
\special{pa 3630 4040}%
\special{pa 3603 4058}%
\special{pa 3573 4070}%
\special{pa 3542 4077}%
\special{pa 3510 4080}%
\special{pa 3478 4079}%
\special{pa 3447 4075}%
\special{pa 3416 4066}%
\special{pa 3387 4053}%
\special{pa 3362 4032}%
\special{pa 3350 4003}%
\special{pa 3350 4000}%
\special{sp}%
%
\special{pn 8}%
\special{ar 3510 4080 560 250  0.0000000 6.2831853}%
%
\special{pn 8}%
\special{pa 4080 330}%
\special{pa 4080 4080}%
\special{fp}%
%
\special{pn 8}%
\special{pa 2950 340}%
\special{pa 2950 4080}%
\special{fp}%
%
\special{pn 8}%
\special{ar 1410 330 560 250  0.0000000 6.2831853}%
%
\special{pn 8}%
\special{ar 1410 4090 560 250  0.0000000 6.2831853}%
%
\special{pn 8}%
\special{pa 1980 340}%
\special{pa 1980 4090}%
\special{fp}%
%
\special{pn 8}%
\special{pa 850 350}%
\special{pa 850 4090}%
\special{fp}%
%
\special{pn 8}%
\special{ar 5890 310 560 250  0.0000000 6.2831853}%
%
\special{pn 8}%
\special{ar 5890 4070 560 250  0.0000000 6.2831853}%
%
\special{pn 8}%
\special{pa 6460 320}%
\special{pa 6460 4070}%
\special{fp}%
%
\special{pn 8}%
\special{pa 5330 330}%
\special{pa 5330 4070}%
\special{fp}%
\put(34.0000,-47.9000){\makebox(0,0)[lb]{Figure 9.3}}%
%
\special{pn 20}%
\special{pa 3350 350}%
\special{pa 3340 2050}%
\special{fp}%
%
\special{pn 20}%
\special{pa 3650 360}%
\special{pa 3650 2050}%
\special{fp}%
%
\special{pn 20}%
\special{pa 3350 3980}%
\special{pa 3350 2390}%
\special{fp}%
%
\special{pn 20}%
\special{pa 3650 3960}%
\special{pa 3650 2400}%
\special{fp}%
%
\special{pn 8}%
\special{pa 3650 2080}%
\special{pa 3650 2390}%
\special{fp}%
%
\special{pn 8}%
\special{pa 3350 2060}%
\special{pa 3360 2400}%
\special{fp}%
%
\special{pn 8}%
\special{ar 3490 2050 150 80  0.0000000 6.2831853}%
\end{picture}%

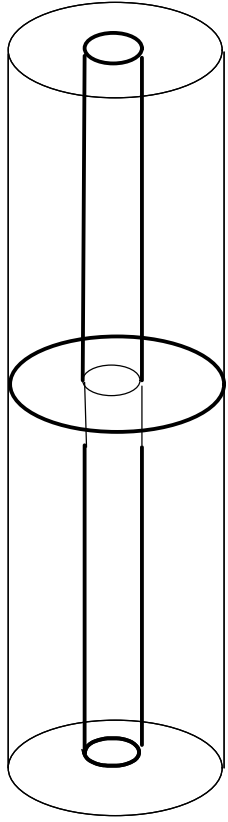
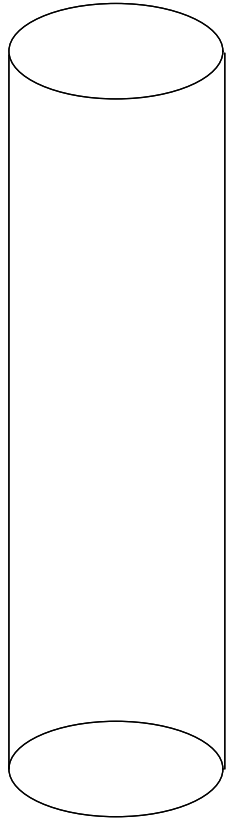
\end{document}